\documentclass[11pt]{article}
\usepackage{psfrag,
graphicx, harvard,
epsf,
euscript,amsfonts,amsmath,latexsym,amssymb}
\newtheorem{lemma}{Lemma}
\newtheorem{theorem}{Theorem}
\newtheorem{proposition}{Proposition}

\newtheorem{remark}{Remark}
\newtheorem{definition}{Definition}

\newcommand{\rE}{\mathbb E}

\newcommand{\rR}{\mathbb R}
\newcommand{\rP}{\mathbb P}

\newcommand{\rC}{\mathbb C}
\newcommand{\rS}{\mathbb S}

\newcommand{\rL}{\mathbb L}

\newcommand{\rF}{\mathbb F}

\newcommand{\rH}{\mathbb H}
\renewcommand{\kappa}{\varkappa}

\newcommand{\qed}{\mbox{\hspace{1em}\rule{2mm}{2mm}}}
\def\noi{\noindent}

\begin{document}
\author{
Alexander Goldenshluger
\thanks{Department of Statistics, University of Haifa, 31905 Haifa, Israel.
 email: {\em goldensh@stat.haifa.ac.il}}   \and 
Oleg Lepski\thanks{Laboratoire d'Analyse, Topologie et Probabilit\'es
UMR CNRS 6632, Universit\'e de Provence,
39, rue F.Joliot Curie,
13453 Marseille, France.
email: {\em lepski@cmi.univ-mrs.fr}}
}
\title{Structural Adaptation via $\rL_p$-norm Oracle Inequalities}

\maketitle

\begin{abstract} 
In this paper we study the problem of adaptive estimation of a multivariate 
function satisfying some structural assumption. We propose 
a novel  estimation procedure 
that adapts simultaneously to unknown structure and  smoothness of  the
underlying function.
The problem of structural adaptation is stated 
as the problem of selection from a given collection of estimators.
We develop a general selection rule and establish for it 
global oracle inequalities under arbitrary $\rL_p$--losses.
These results are applied for adaptive estimation 
in the additive multi--index model.  
\end{abstract}

\vspace*{1em} \noi {\bf Short Title:} Structural adaptation via oracle 
inequalities

\noi {\bf Keywords:} structural adaptation, oracle inequalities, minimax risk,
adaptive estimation, optimal rates of convergence

\noi {\bf 2000 AMS Subject Classification} : 62G05, 62G20

%\fi

%%%%%%%%%%%%%%%%%%%%
%\setlength{\evensidemargin}{-0.5cm}
%\setlength{\oddsidemargin}{-0.5cm} \setlength{\topmargin}{-1cm}
\def\huh{\hbox{\vrule width 2apt height 8apt depth 2apt}}
\def\Bbb#1{{\bf #1}}
\def\eqnum#1{\eqno (#1)}
\def\fnote#1{\footnote}
\def\blacksquare{\hbox{\vrule width 4apt height 4apt depth 0apt}}
\def\square{\hbox{\vrule\vbox{\hrule\phantom{o}\hrule}\vrule}}
\def\inter{\mathop{{\rm int}}}
\def\epi{{\rm epi}}
\def\nr{\par \noindent}
\def\diag{{\rm diag}}
\def\trace{{\rm Trace}}
\def\clsr{{\rm cl}}
\def\beq{\begin{equation}}
\def\eeq{\end{equation}}
\newcommand{\ha}{\widehat{\alpha}}
\newcommand{\hbt}{\widehat{\beta}}
\newcommand{\tbt}{\tilde{\beta}}
\newcommand{\htf}{\widehat{f}}
\newcommand{\ljk}{\Lambda_{jk}}
\newcommand{\tj}{\Theta_j}
\newcommand{\downsamp}{DWNS}
% Greek letters 
\newcommand{\al}{\alpha}
\newcommand{\bt}{\beta}
\newcommand{\g} {\gamma}
\newcommand{\wh}[1] {\widehat{#1}}
\newcommand{\gt}{\frac{\gamma}{t}}
\newcommand{\de}{\overline{\delta}}
\newcommand{\dl}{\delta}
\newcommand{\sg}{\sigma}
\newcommand{\z}{\zeta}
\newcommand{\bD}{\bar{\Delta}}
\newcommand{\bP}{\bar{\Pi}}
\newcommand{\D}{\Delta}
\newcommand{\bd}{\bar{\Delta}}
\newcommand{\bw}{\bar{w}}
\newcommand{\bG}{\bar{G}}
\newcommand{\te}{\theta}
\newcommand{\hte}{\hat{\theta}}
\newcommand{\vth}{\vartheta}
\newcommand{\vf}{\varphi}
\newcommand{\e}{\varepsilon}
\newcommand{\ups}{\upsilon}
\newcommand{\Up}{\Upsilon}
\newcommand{\ve}{\varepsilon}
\newcommand{\we}{\wedge}
\newcommand{\ot} {\frac{1}{t}}
\newcommand{\on} {\frac{1}{n}}
\newcommand{\rl}[1]{~\cite{#1}}
\newcommand{\G} {\bar{G}}
\newcommand{\bpsi}{\mbox{\boldmath $\psi$}}
\newcommand{\bphi}{\mbox{\boldmath $\phi$}}
\newcommand{\supp}{\mbox{supp}}
\newcommand{\diamm}{\mbox{diam}}
\newcommand{\spann}{\mbox{span}}
\newcommand{\sign}{\mbox{sign}}
\newcommand{\lan}{\langle}
\newcommand{\ran}{\rangle}
\newcommand{\tm}{\times}
\newcommand{\mt}{m(T)}
\newcommand{\nb}{\nabla}
\newcommand{\ncc}{\nabla^{2}}
\newcommand{\lmt}{\lim_{t \rightarrow \infty}}
\newcommand{\lmn}{\lim_{n \rightarrow \infty}}
\newcommand{\llm} {\lim_{\overline{t \rightarrow \infty}}}
\newcommand{\rt}{\rightarrow}
\newcommand{\uli}[1]{\overline{\lim_{#1}}}
\newcommand{\cP}{\stackrel{P}{\rightarrow}}
\newcommand{\defn}{\stackrel{\Delta}{=}}
\newcommand{\Ey}{E_{y,a}}
\newcommand{\Eb}{\bar{E}_{U,n}}
\newcommand{\cD}{{\cal D}}
\newcommand{\cPP}{{\cal P}}
\newcommand{\cO}{{\cal O}}
\newcommand{\cF}{{\cal F}}
\newcommand{\cC}{{\cal C}}
\newcommand{\cE}{{\cal E}}
\newcommand{\cK}{{\cal K}}
\newcommand{\cG}{{\cal G}}
\newcommand{\cN}{{\cal N}}
\newcommand{\cA}{{\cal A}}
\newcommand{\cU}{{\cal U}}
\newcommand{\cI}{{\cal I}}
\newcommand{\cp}{{\cal P}}
\newcommand{\cB}{{\cal B}}
\newcommand{\cM}{{\cal M}}
\newcommand{\Pm}{\hbox{\it I\hskip -3apt P}}
\newcommand{\Em}{\hbox{\it I\hskip -3apt E}}
\newcommand{\mI}{{\bf I}}
\newcommand{\Cm}{{\bf C}}
\newcommand{\Rm}{{\bf R}}
\newcommand{\Nm}{{\bf N}}
\newcommand{\Zm}{{\bf Z}}
\newcommand{\be}{\begin{eqnarray}}
\newcommand{\ee}[1]{\label{eq:#1}\end{eqnarray}}
\newcommand{\nn}{\nonumber \\}
\newcommand{\ese}{\end{eqnarray*}}
\newcommand{\bse}{\begin{eqnarray*}}
\newcommand{\rf}[1]{~(\ref{eq:#1})}
\newcommand{\sm}[2]{\sum_{i=#1}^{#2}}
\newtheorem{thee}{Theorem}
\newcommand{\bth}{\begin{thee}}
\newcommand{\rth}[1]{\ref{the:#1}}
\newtheorem{col}{Corollary}
\newcommand{\bcol}{\begin{col}}
\newcommand{\ecol}[1]{\label{the:#1}\end{col}}
\newcommand{\rc}[1]{\ref{the:#1}}
\newtheorem{defi}{Definition}
\newcommand{\bdf}{\begin{defi}}
\newcommand{\edf}[1]{\label{df:#1}\end{defi}}
\newcommand{\rdf}[1]{\ref{df:#1}}
\newtheorem{lem}{Lemma}
\newcommand{\blem}{\begin{lem}}
\newcommand{\elem}[1]{\label{le:#1}\end{lem}}
\newcommand{\rle}[1]{\ref{le:#1}}
\newtheorem{pro}{Proposition}
\newcommand{\bpro}{\begin{pro}}
\newcommand{\epro}[1]{\label{pr:#1}\end{pro}}
\newcommand{\rp}[1]{\ref{pr:#1}}
\newtheorem{exmp}{Example}
\newcommand{\bexmp}[1]{\begin{exmp}\label{exmp:#1}\end{exmp}}
\newcommand{\eexmp}[1]{\label{exmp:#1}\end{exmp}}
\newcommand{\rexmp}[1]{\ref{exmp:#1}}
\newtheorem{exos}{Exercise}
\newcommand{\bexo}[1]{\begin{exos}\label{exo:#1}\end{exos}}
%\newcounter{exose}[chapter]
%\newcommand{\bexo}[1]{\refstepcounter{exose}[chapter]\label{exo:#1}\paragraph{Exercice \arabic{chapter}.\arabic{exose}}}
%\newcommand{\rexo}[1]{\ref{exo:#1}}
%
%\newcounter{assc}
%\newcommand{\bass}[1]{\refstepcounter{assc}\label{ass:#1} \begin{maliste}{\bf \arabic{assc}.}}
%\newcommand{\eass}{\end{maliste}}
%\newcommand{\ras}[1]{\ref{ass:#1}}
%
\newcommand{\pr} {\par \noindent{\bf Proof\,:~}}
\newcommand{\epr}{\hfill\hbox{\hskip 4pt
                \vrule width 5pt height 6pt depth 1.5pt}\vspace{0.5cm}\par}
\def    \qed    {\hfill\hbox{\hskip 4pt
                \vrule width 3pt height 1pt depth 1pt
                \hbox{\vrule width 2pt height 3.5pt depth 1pt}}}
\newcommand{\pht}{\hat{\phi}}

\newcommand{\om}{\omega}
\newcommand{\wf}{\hat{\phi}}
\newcommand{\whf}{\hat{f}}
\newcommand{\wbb}{\bar{\beta}}
\newcommand{\bal}{\bar{\alpha}}
\newcommand{\tal}{\tilde{\alpha}}
\newcommand{\tbe}{\tilde{\beta}}
\newcommand{\hal}{\hat{\alpha}}
\newcommand{\card}{\mbox{card}}
\newcommand{\wA}{\hat{A}}
\newcommand{\osc}{\mbox{osc}}
\newcommand{\sR}{{\sf R}}
%
%\def\baselinestretch{2.0}
%
% Arik's macros
\def\huh{\hbox{\vrule width 2pt height 8pt depth 2pt}}
\def\Bbb#1{{\bf #1}}
\def\eqnum#1{\eqno (#1)}
\def\fnote#1{\footnote}
\def\pn{\par\noindent}
\def\blacksquare{\hbox{\vrule width 4pt height 4pt depth 0pt}}
\def\square{\hbox{\vrule\vbox{\hrule\phantom{o}\hrule}\vrule}}
\def\inter{\mathop{{\rm int}}}
\def\epi{{\rm epi}}
\def\nr{\par \noindent}
\def\diag{{\rm diag}}
\def\trace{{\rm Trace}}
\def\clsr{{\rm cl}}
\def\beq{\begin{equation}}
\def\eeq{\end{equation}}
\def\eqbd{\colon =}
\renewcommand\arraystretch{1}
\def\ba{\begin{array}}
\def\ea{\end{array}}
\def\Var{\mbox{Var}}
\def\Cov{\mbox{Cov}}
\def\Corr{\mbox{Corr}}
\def\beann{\begin{eqnarray*}}
\def\eeann{\end{eqnarray*}}
\def\bea{\begin{eqnarray}}
\def\eea{\end{eqnarray}}
\def\Tr{\mathop{{\rm Tr}}}
\def\Det{\mathop{{\rm Det}}}
\def\Sz{\mathop{{\hbox{\sl Size}}}}
\def\Vl{\mathop{{\hbox{\sl Vol}}}}
\def\CE{\mathop{{\hbox{\sl CircEll}}}}
\def\CS{\mathop{{\hbox{\sl CircSimp}}}}
\def\IE{\mathop{{\hbox{\sl InscrEll}}}}
\def\nrm{\mathop{{\parallel}}}
\def\dim{\mathop{{\rm dim}\,}}
\def\inter{\mathop{{\rm int}\,}}
\def\vrt{\mathop{\,\mid}\,}
\def\argmin{\mathop{\rm argmin}}
\def\argmax{\mathop{\rm argmax}}
\def\cl{\mathop{{\rm cl}}}
\def\Dom{\mathop{{\rm Dom}}}
\def\Ker{\mathop{{\rm Ker}}}
\def\conv{\mathop{{\rm Conv}}}
\def\Diag{\mathop{{\rm Diag}}}
\def\Rang{\mathop{{\rm Rang}}}
\def\Compl{\mathop{{\rm Compl}}}
\def\Accur{\mathop{{\rm Accur}}}
\def\log{\mathop{{\rm log}}}
\def\mes{\mathop{{\rm mes}}}
\def\Vol{\mathop{{\rm Vol}}}
\def\vol{\mathop{{\rm vol}}}
\def\dist{\mathop{{\rm dist}}}
\def\Argmin{\mathop{\rm Argmin}}
\def\size{\mathop{{\rm Size}}}
\def\Diam{\mathop{{\rm Diam}}}
\def\AvDiam{\mathop{{\rm AvDiam}}}
\def\Residual{\mathop{{\rm Residual}}}
\def\Size{\mathop{{\rm Size}}}
\def\me{\mathop{{\rm e}}}
\def\sign{\mathop{{\rm sign}}}
\def\EllOut{\mathop{{\rm EllOut}}}
\def\Ell{\mathop{{\sl Ell}}}
\def\card{\mathop{{\rm card}}}
\def\book{\cite{NNbook}}
\def\Kappa{{\mathop{\,{\bar\kappa}}}}
\def\Pol{{\mathop{\rm Polyn}}}
\def\Seq{{\mathop{\rm Seq}}}
%%%%
% bold math
%%%%
\newcommand{\blh}{\boldsymbol{h}}
\newcommand{\blx}{\boldsymbol{x}}
\newcommand{\bly}{\boldsymbol{y}}
\newcommand{\blu}{\boldsymbol{u}}
\newcommand{\blz}{\boldsymbol{z}}
\newcommand{\blt}{\boldsymbol{t}}
\newcommand{\blk}{\boldsymbol{k}}
\newcommand{\blU}{\boldsymbol{U}}
\newcommand{\bltheta}{\boldsymbol{\theta}}
\newcommand{\blTheta}{\boldsymbol{\Theta}}
\newcommand{\bltau}{\boldsymbol{\tau}}
\newcommand{\blJ}{\boldsymbol{J}}

\renewcommand{\overline}{\bar}
%\renewcommand{\blh}{h}
%\renewcommand{\blx}{x}
%\renewcommand{\bly}{y}
%\renewcommand{\blu}{u}
%\renewcommand{\blz}{z}
%\renewcommand{\blt}{t}
%\renewcommand{\blk}{k}
%\renewcommand{\bltheta}{\theta}
%\renewcommand{\blTheta}{\Theta}
%\renewcommand{\bltau}{\tau}
%\renewcommand{\blJ}{J}

%
% Macros for the TDs
%
\def\vi{\stackrel{\to}{i}}
\def\vj{\stackrel{\to}{j}}
\def\vv{\stackrel{\to}{v}}
\def\va{\stackrel{\to}{a}}
\def\vk{\stackrel{\to}{k}}
\def\vb{\stackrel{\to}{b}}
\def\vc{\stackrel{\to}{c}}
\def\vn{\stackrel{\to}{n}}
\def\vd{\stackrel{\to}{d}}
\def\vg{\stackrel{\to}{g}}
\def\vr{\stackrel{\to}{r}}
\def\vom{\stackrel{\longrightarrow}{OM}}
\def\vgrad{\stackrel{\longrightarrow}{\mbox{grad}}}
\def\vu{\stackrel{\to}{u}}
\def\vF{\stackrel{\to}{F}}
\def\d{\partial}
\def\dx#1{{\partial #1\over \partial x}}
\def\dr#1{{\partial #1\over \partial r}}
\def\dte#1{{\partial #1\over \partial \theta}}
\def\dphi#1{{\partial #1\over \partial \phi}}
\def\drho#1{{\partial #1\over \partial \rho}}
\def\dxx#1{{\partial^2 #1\over \partial x^2}}
\def\dxy#1{{\partial^2 #1\over \partial x \partial y}}
\def\dyx#1{{\partial^2 #1\over \partial y \partial x}}
\def\dy#1{{\partial #1\over \partial y}}
\def\dyy#1{{\partial^2 #1\over \partial y^2}}
\def\dyz#1{{\partial^2 #1\over \partial y\partial z}}
\def\dxz#1{{\partial^2 #1\over \partial x\partial z}}
\def\dz#1{{\partial #1\over \partial z}}
\def\dzz#1{{\partial^2 #1\over \partial z^2}}
\def\grad{\stackrel{\longrightarrow}{\mbox{grad}}}
\def\indep{\hbox{$\perp\hskip -4apt \perp$}}
\def\qr{\qquad\cr}
%\def\figure{$$\begin{array}{c} \qr \qr \qr \qr \qr \qr \qr \qr \qr \qr \qr  \end{array}$$}
%%%%%%%%%%%%%%%%%%%%%%%%%%%%%%%%%%%%%%5
%\begin{document}
%%%%%%%%%%%%%%%%%%%%%%%%%%%%%%%%%%%%%%
%\renewcommand{\baselinestretch}{1.3}  % increases spacing x1.3
%---------------------------
\section{Introduction}
\subsection{Motivation}
In this paper we study the problem of minimax 
adaptive estimation of an unknown function $F:\rR^d\to\rR$ in the
multidimensional Gaussian white noise model
\begin{equation}\label{model} Y(dt)=F(t)dt +\e W(dt), \quad
t=(t_1,\ldots,t_d) \in {\cal D} ,
\end{equation}
 where $\cD\supset [-1/2,1/2]^d$ is an open interval in
$\rR^d$, $W$ is the standard Brownian sheet in $\rR^d$ and
$0<\e<1$ is the noise level. Our goal is to estimate the function
$F$ on the set ${\cal D}_0:=[-1/2,1/2]^d$ 
from the observation $\{Y(t), \ t\in
\cD\}$.
 We consider the observation set $\cD$ which is larger than $\cD_0$ in
order to avoid discussion of boundary effects.
We would like to emphasize that
such assumptions are rather usual in multivariate models, see, e.g.,
\citeasnoun{Hall} and \citeasnoun{Chen}.
\par
To measure performance of estimators, we will use the risk
function determined by the $\rL_p$-norm $\|\cdot \|_p, 1\leq p\leq\infty$ on
$\cD_0$: for $F:\rR^d\to\rR$, $0<\varepsilon<1$,  and for
an arbitrary estimator $\tilde{F}$ based on the
observation $\{Y(t), \ t\in \cD\}$ we consider the risk
$$
{\cal R}_p [\tilde{F}; F] =
\rE_F||\tilde{F}-F||_p.
$$
Here and in what follows $\rE_F$ denotes the expectation with respect to the
distribution $\rP_F$ of the observation $\{Y(t), \ t\in \cD\}$
satisfying (\ref{model}).
\par
We will suppose that $F\in\cG_s$, where $\{\cG_s,
s\in{\cal S}\}$ is a collection of functional classes indexed by
$s\in {\cal S}$. The choice of this collection is a delicate problem,
and below we discuss it in detail.
\par
For a given class $\cG_s$ we define the maximal risk
\begin{equation}
\label{mmrisk}
{\cal R}_p[\tilde{F}; \cG_s]=
\sup_{g\in\cG_s} {\cal R}_p [\tilde{F};F],
\end{equation}
and 
% we 
study  asymptotics (as the noise level
$\varepsilon$ tends to 0) of the minimax risk
\begin{equation*}
%\label{mmrisk}
\inf_{\tilde{F}}
{\cal R}_p [\tilde{F};\cG_s]
\end{equation*}
where $\inf_{\tilde{F}}$ denotes the infinum over all
estimators of $F$. At this stage, we suppose that parameter $s$ is
known, and therefore the functional class $\cG_s$ is fixed. In
other words, we are interested in minimax estimation of $F$. 
The important remark in this context is that
the minimax rate of convergence $\phi_{\varepsilon}(s)$ on
$\cG_s$ (the rate which satisfies
$\phi_{\varepsilon}(s)\asymp
\inf_{\tilde{F}}{\cal R}_p[\tilde{F};\cG_s]$)
as well as  the estimator attaining this rate  (called the rate
optimal estimator in asymptotic minimax sense) depend
on parameter $s$. This dependence  restricts 
application of 
the minimax approach in practice.
Therefore, our main goal is to construct an estimator which is
independent of $s$ and achieves the minimax rate
$\phi_{\varepsilon}(s)$ simultaneously for all
$s\in {\cal S}$. Such an estimator, if it exists, is called optimally
adaptive on ${\cal S}$. 
\par
Let us discuss now the choice of the collection 
$\{\cG_s, s\in {\cal S}\}$.
It is well known that the main difficulty in estimation of
multivariate functions is the curse of dimensionality:
the best
attainable rate of convergence of  estimators becomes very
slow, as the dimensionality grows.
To illustrate this effect, suppose,
for example, that the underlying function $F$ belongs to
$\cG_s=\rH_d(\alpha,L)$, $s=(\alpha, L)$, $\alpha>0, L>0$, where
$\rH_d(\alpha,L)$ is an isotropic H\"older ball of functions. We
give the exact definition of this functional class later. Here we
only mention that $\rH_d(\alpha,L)$ consists of functions $g$ with
bounded partial derivatives of order $\leq\lfloor\alpha\rfloor$ and
such that, for all $x,y\in\cD$,
$$
\left|g(y)-P_{g}(x,y-x)\right|\leq L|x-y|^\alpha,
$$
where $P_{g}(x,y-x)$ is the Taylor polynomial of order
$\leq\lfloor\alpha\rfloor$ obtained by expansion of $g$ around the
point $x$, and $|\cdot |$ is the Euclidean norm in $\rR^d$.
Parameter $\alpha$ characterizes the isotropic (i.e., the same in
each direction) smoothness of function $g$.
\par
If we use the risk (\ref{mmrisk}), uniformly on $\rH_d(\alpha,L)$
the rate of convergence of estimators cannot be asymptotically
better than
 \begin{equation}\label{eq:psi}
\psi_{\varepsilon,d}(\alpha)=\left\{\begin{array}{ll} 
\e^{2\alpha/(2\alpha+d)}, & p\in [1, \infty)\\
(\e\sqrt{\ln\e^{-1}})^{2\alpha/(2\alpha+d)} & p=\infty.
                                    \end{array}
\right.
\end{equation}
[cf. \citeasnoun{ibr-has82}, \citeasnoun{stone82}, 
\citeasnoun{nussbaum87},  \citeasnoun{Bertin}]. 
This is the minimax rate on
$\rH_d(\alpha,L)$: in fact, it can be achieved by a kernel
estimator with properly chosen bandwidth and kernel. More general
results on asymptotics of the minimax risks in estimation of
multivariate functions can be found in 
\citeasnoun{lepski-picard} and
\citeasnoun{Bertin}. 
It is clear that if $\alpha$ is fixed 
then even for moderate $d$ the estimation accuracy 
is very poor unless the noise level $\e$ is unreasonably small.
%and $d$ is large enough this
%asymptotics 
%of the rate 
%is too pessimistic to be used for real data:
%the value $\psi_{\varepsilon,d}(\alpha)$ is small only if the noise
%level $\varepsilon$ is unreasonably small. On the other hand, if the
%noise level $\varepsilon$ is realistically small the above
%asymptotics might be of no use already in dimension $2$ or $3$.
\par
This problem arises because the $d$-dimensional H\"older ball
$\rH_d(\alpha,L)$ is too massive.  A way to overcome the curse of
dimensionality is to consider models with 
%poorer 
smaller functional classes
$\cG_s$. Clearly, if the class of candidate functions $F$ is
smaller, the rate of convergence of estimators is faster. Note that
the ``poverty'' of a functional class can be described in terms of
restrictions on its metric entropy. There are 
nevertheless 
several ways to do it.
\subsection{Structural adaptation}
In this paper we will follow the modeling strategy which consists in imposing
 additional structural assumptions on the function 
to be estimated. This approach was pioneered by \citeasnoun{stone85}
who discussed the trade--off between flexibility and dimensionality of
nonparametric models
and formulated {\em the heuristic dimensionality reduction principle}.
The main idea is to assume that even though $F$ is a $d$--dimensional
function, it has a simple structure such that $F$
is effectively $m$--dimensional with $m<d$. 
The standard examples
of structural nonparametric models are the following.
\begin{itemize}
\item[(i)]  {\em [Single--index model.]} Let $e$ be a direction vector
in $\rR^d$, and assume that $F(x)=f(e^Tx)$ for some 
unknown univariate function $f$. 
\item[(ii)] {\em [Additive model.]}. Assume that
$F(x)=\sum_{i=1}^d f_i(x_i)$, where  
$f_i$ are unknown univariate functions.
\item[(iii)] {\em [Projection pursuit regression.]} Let $e_1, \ldots, e_d$
be direction vectors in $\rR^d$, and assume that
$F(x)=\sum_{i=1}^d f_i (e^Tx)$, where $f_i$ are as in (ii).
\item[(iv)] {\em [Multi--index model.]} Let $e_1, \ldots, e_m$, $m<d$ 
are direction vectors
and assume that $F(x)=f(e_1^Tx, \ldots, e_m^T x)$ for some unknown $m$-dimensional
function $f$.
\end{itemize}
In the first three examples the function $F$ is effectively one--dimensional,
while in the fourth one it is $m$--dimensional.
The heuristic dimensionality reduction principle by \citeasnoun{stone85} 
suggests that the optimal rate of convergence attainable in 
structural nonparametric models should
correspond to the effective dimensionality of $F$.
\par
Let us make the following important remark.
% which will be useful in the sequel.
\par
 The estimation problem in the models of types (i), (iii) and (iv)
 can be viewed as the problem of  adaptation to unknown 
 structure (structural adaptation). 
Indeed, 
%let us suppose that 
if the direction vectors are given
 then, 
after 
%simple 
a linear transformation, 
the problem is 
%can be 
reduced either to the estimation problem in
 the additive model (cases~(i) and (iii)) 
or to the estimation of an $m$-variate function. This explains 
 the form of minimax rate of convergence.
The main problem however is to find an estimator that
% would 
adjusts automatically to unknown
 direction vectors. 
\par
This remark allows to state 
the problem of structural adaptation in the following 
rather general way.
\subsection{$\rL_p$--norm oracle inequalities}
Suppose that we are given a collection of estimators
 $\big\{\hat{F}_{\theta}, \theta\in\Theta\subset\rR^m \big\}$
 based on the observation  $\{Y(t), \ t\in\cD\}$. 
 In the previous examples parameter $\theta$ could be, for instance,
 the unknown matrix $E=(e_1,\ldots,e_d)$ of the direction vectors, $\theta=E$,
 and $\hat{F}_{E}$ could be a kernel estimator 
 constructed under hypothesis that $E$ and smoothness of 
the functional components are known 
(a kernel estimator with fixed bandwidth).
\par 
% To 
With each estimator $\hat{F}_\theta$ and 
%the 
unknown function $F$
we associate the risk 
 ${\cal R}_p[\hat{F}_{\theta};F]$. The problem 
%we address now consists in 
is
% constructing 
to construct an estimator, say, $\hat{F}_*$ such that 
for all $F$ obeying 
 given smoothness conditions one has
 \begin{equation}
 \label{oracle}
 {\cal R}_p [\hat{F}_*;F] \leq {\cal L}\inf_{\theta\in\Theta}
 {\cal R}_p [\hat{F}_{\theta};F],
 \end{equation}
where $\cal L$ is an absolute constant independent of $F$ and $\e$. 
Following the modern statistical 
% language 
terminology we will call the 
inequality (\ref{oracle}) 
{\it the $\rL_p$-norm oracle inequality}.
\par
Returning to our example with $\theta=E$ we observe that being established,  
the {\it $L_p$-norm oracle inequality} leads immediately to the minimax result
for any given value of smoothness parameter $(\alpha,L)$.
In particular, we can state that the estimator $\hat{F}_*$ is 
adaptive with respect to unknown structure.
\par
It is important to realize that
the same strategy 
%permits  
allows to avoid 
dependence of estimation procedures on 
%the 
smoothness. 
%For that 
To this end it is sufficient
 \begin{itemize}
 \item
 to consider $\theta=(E,\alpha,L)$ 
that leads to the collection of kernels estimators
 with the non-fixed bandwidth and orientation;
 \item
  to propose 
an estimator $\hat{F}_*$ based on this collection;
  \item
  to establish for this estimator 
{\it the $L_p$-norm oracle inequality} (\ref{oracle}) for
  any $F\in L_2(\cD)$ (or on a bit 
%narrow 
smaller functional space).
 \end{itemize}
Being realized, 
this program leads to 
%the 
an  estimator 
% which 
that is adaptive with respect to
unknown structure 
%as well as to 
and unknown smoothness properties. 
It is important to note 
that such methods allow  to estimate 
multivariate functions with high accuracy 
without sacrificing flexibility of 
%the 
modeling.
\subsection{Objective of the paper}
The goal of the present paper is at least two--fold. 
\par
First we introduce
and study a general structural model that we  call {\em the additive 
multi--index model}; it includes models (i)--(iv) 
as special cases.
This generalization is dictated by the following reasons. On the one hand,
structural assumptions allow to improve 
the quality of statistical analysis. On the other hand, they can lead to inadequate
modeling. Thus we seek  a general structural model that still
allows to gain in estimation accuracy.  To our knowledge the additive
multi--index model
did not previously appear in the statistical literature.
For this model we propose an estimation procedure that adapts
simultaneously to unknown structure and smoothness of the 
underlying function.
% (structural adaptation). 
The adaptive results are
obtained for $\rL_\infty$--losses and for a scale  of the 
H\"older type functional classes.
\par
To study this model we proceed as follows. We state the problem of structural 
adaptation as the problem of 
selection from a given collection of estimators.
For a collection of linear estimators satisfying rather mild assumptions
we propose a novel general selection rule and establish for it
the $\rL_p$--norm oracle inequality (\ref{oracle}). 
Similar ideas were used 
in \citeasnoun{lepski-levit},
\citeasnoun{lepski-picard}, \citeasnoun{lepski-composite}
for {\em pointwise} adaptation.
However we emphasize that
our work is the first where  
the $\rL_p$--norm oracle inequality
is derived directly
without applying pointwise estimation results.  
%without passing to pointwise estimation.
It is precisely this fact that
allows to obtain adaptive results for 
arbitrary $\rL_p$--losses.
The selection rule as well as the $\rL_p$--norm oracle inequality
are not related 
to any specific model, and they 
are applicable 
in a variety of setups where
linear estimators are appropriate. 
We apply these general 
results to a specific collection of kernel estimators corresponding to the
additive multi--index model.
\subsection{Connection to other works}
\paragraph*{\underline{Structural models.}}
The heuristic dimensionality
reduction principle was proved by 
\citeasnoun{stone85} for the additive model~(ii),
and by \citeasnoun{Chen} 
and \citeasnoun{golubev}
for the projection pursuit regression model~(iii) that 
includes as a particular case the single--index model~(i).
In particular, it was shown there that in these models 
the asymptotics of the risk
(\ref{mmrisk}) with $p=2$ and with $\cG_{s}$, $s=(\alpha,L)$, where   $\cG_{s}$ 
is either H\"older or 
Sobolev ball,
is given by $\psi_{\varepsilon,1}(\alpha)$. As we see, the accuracy of estimation in
such models   corresponds to the one--dimensional rate ($d=1$).
\par
Further results and references on estimation 
in models (i)--(iv) can be found, e.g., in 
\citeasnoun{yatracos},
\citeasnoun[Chapter~22]{gyorfi}, and
\citeasnoun{ibragimov2004}.
Let us briefly discuss the results obtained.
\begin{itemize}
\item
The  estimators providing  
the rate mentioned above depend heavily on the use of 
 $\rL_2$-losses ($p=2$) 
%for the description of 
in the risk definition. 
As a consequence, all proposed constructions
 cannot be used for any other types of loss functions.
 \item
 Except for the paper by 
\citeasnoun{golubev}, where the estimator independent on the parameter
  $s=(\alpha,L)$ was proposed for the 
model~(i), all other estimators depend explicitly 
  on the prior information on smoothness of the underlying function.
\item 
 As far as we know there are no even 
minimax results obtained for the model~(iv). One can guess that
 asymptotics of the risk (\ref{mmrisk}) is given 
by $\psi_{\varepsilon, m}(\alpha)$ which is much better then
 $d$-dimensional rate $\psi_{\varepsilon,d}(\alpha)$ since $m<d$.
\end{itemize}
It is also worth mentioning that
there is vast literature on estimation of vectors $e_i$,  when $f_i$
are treated as nonparametric nuisance parameters;
see, e.g.,
\citeasnoun{huber}, \citeasnoun{Hall}, 
\citeasnoun{juditsky1} and references therein.
\paragraph*{\underline{Oracle approach.}}
To understand the place of the oracle approach 
within the theory of nonparametric estimation
let us quote
\citeasnoun{johnstone}: 
\begin{quote}
{\em
"Oracle inequalities are neither the beginning 
nor the end of a theory, but when available, are informative tools." 
}
\end{quote}
Indeed, oracle inequalities are 
% at least  
very powerful tools for deriving
minimax and minimax adaptive results.
The aim of the oracle approach can be formulated
as follows: given a collection of different estimators based on available
data, select the best estimator from the family
({\em model selection})
[see, e.g., \citeasnoun{barron-birge-massart}], 
or find the best convex/linear combination
of the estimators from the  family ({\em convex/linear aggregation})
[see \citeasnoun{nemirovski2000}, \citeasnoun{tsybakov2003}].
The formal  definition of the oracle 
requires specification of
the collection of estimators and
the criterion of optimality.
\par
The majority of oracle procedures described in the literature
use the $\rL_2$--risk
as the criterion of optimality.
The following methods can be cited in this context:
penalized likelihood estimators,
unbiased risk estimators, blockwise Stein estimators, risk hull
estimators  and so on
[see  \citeasnoun{barron-birge-massart}, \citeasnoun{cavalier},
\citeasnoun{golubev-2004}
and references therein].
The most general results in the framework of
$\rL_2$--risk aggregation theory were obtained by
\citeasnoun{nemirovski2000}
who showed how to aggregate 
arbitrary estimators. 
%Let us especially mention the results 
%on $\rL_2$ model selection obtained by Birg\'e, Massart and
%their school. 
\par
Other  oracle procedures were developed in the context of 
pointwise estimation; see, e.g., 
\citeasnoun{lepski-spokoiny}, \citeasnoun{gnem},
\citeasnoun{spok-belom} 
for the univariate case, and
 \citeasnoun{lepski-levit}, 
\citeasnoun{lepski-picard} for the multivariate case.
Moreover, 
%in the papers 
\citeasnoun{lms} and
\citeasnoun{lepski-picard} 
%it was shown 
show how to derive $\rL_p$--norm oracle inequalities from pointwise
oracle inequalities. Although
these $\rL_p$--norm oracle
inequalities allow to derive minimax results on rather complicated
functional spaces, they do not lead to sharp adaptive results.
\par
Finally we mention the $\rL_1$--norm oracle approach developed
by \citeasnoun{devroye-lugosi} in context of density estimation.
\par
The rest of the paper is organized as follows. 
In Section~\ref{sec:scheme} we present our general selection rule 
and establish 
%for it 
the {\em key oracle inequality}. 
Section~\ref{sec:multiindex}
is devoted to adaptive estimation in the additive multi--index model.
The proofs of the mains results are given in Section~\ref{sec:proofs}.
Auxiliary results are postponed to Appendix. 

\section{General selection rule}\label{sec:scheme}
\subsection{Preliminaries}
\par
In what follows
$\|\cdot\|_p$ stands for the $\rL_p({\cal D})$--norm, while
$\|\cdot\|_{p,q}$   denotes the
$\rL_{p,q}({\cal D}\times {\cal D}_0)$--norm:
\[
\|G\|_{p,q} = 
\Bigl(\int \Bigl(\int |G (t, x)|^p dt\Bigr)^{q/p} dx\Bigr)^{1/q},\;\;\;
p, q\in [1, \infty].
\] 
We write also $|\cdot|$ for the Euclidean norm,
and it will be always
clear from the context which Euclidean space is meant.
\par
Let $\Theta\subset \rR^m$.
Assume that we are given a 
parameterized family
of kernels 
${\cal K}=\{K_\theta(\cdot, \cdot), \theta\in \Theta\}$, where
$K_\theta : {\cal D}\times {\cal D}_0\to \rR$.
Consider the collection of linear estimators of $F$ associated
with family ${\cal K}$:
\[
{\cal F}({\cal K})=\Bigl\{\hat{F}_\theta (x) = \int K_\theta (t, x) Y(dt),\;\theta\in\Theta\Bigr\}.
\]
Our goal is to propose a measurable choice from the collection 
$\{\hat{F}_\theta, \theta\in \Theta\}$  such that
the risk of the selected estimator will be as close as possible to 
$\inf_{\theta\in \Theta}{\cal R}_p[\hat{F}_\theta, F]$.

\par
Let
\begin{equation}\label{eq:Z-theta}
B_\theta(x) := \int K_\theta (t, x) F(t) dt - F(x),\;\;\;
Z_\theta (x) := \int K_\theta (t, x) W(dt);
\end{equation}
then $\hat{F}_\theta (x) - F(x) = B_\theta(x) +\e Z_{\theta}(x)$, so that
$B_{\theta}(\cdot)$ and $\e Z_{\theta}(\cdot)$ are 
the bias and the stochastic error of the estimator $\hat{F}_\theta$ 
respectively.
We assume that the family ${\cal K}$ of kernels satisfies the following 
conditions.
\begin{itemize}
\item[{\bf (K0)}]
{\em 
For every $x\in {\cal D}_0$ and $\theta\in \Theta$ 
the support of 
$K_\theta (\cdot, x)$  belongs to ${\cal D}$, 
\begin{eqnarray}
\int K_\theta (t, x) dt &=& 1,\;\;\;\forall (x, \theta) \in 
{\cal D}_0\times\Theta,
\label{eq:K1}
\\
\sigma({\cal K})&:=& \sup_{\theta\in \Theta}\|K_{\theta}\|_{2, \infty}
  <\infty,
\label{eq:sigma}
\\
M({\cal K}) &:=& \sup_{\theta\in \Theta}
\Big\{ \sup_x\|K_\theta(\cdot, x)\|_1\; \vee\; \sup_t\|K_\theta (t, \cdot)\|_1
\Big\} < \infty.
\label{eq:M}
\end{eqnarray}
}
\end{itemize}
\begin{remark}
Conditions (\ref{eq:K1}) and (\ref{eq:sigma}) are absolutely
standard in the context of kernel estimation, and only
condition (\ref{eq:M}) has to be discussed. 
First we note that (\ref{eq:M}) is rather mild. 
In particular, if collection ${\cal K}$
contains positive kernels then $M({\cal K})=1$. Moreover 
$M({\cal K})$ will appear in the expression of 
the constant ${\cal L}$ in the {\em $\rL_p$--norm
oracle inequality} (\ref{oracle}).
\end{remark}
\begin{itemize}
\item[{\bf (K1)}]
{\em
For any $\theta, \nu\in \Theta$
\begin{equation}\label{eq:K2}
\int K_\theta (t, y) K_\nu (y, x)dy =
\int K_\theta (y, x) K_\nu (t, y) dy,\;\;\;\forall (x, t)\in {\cal D}_0\times {\cal D}.
\end{equation}
}
\end{itemize}
\begin{remark}
Assumption K1 
is crucial for the construction of our estimation procedure, and 
it restricts the collection of kernels to be used.  
We note nevertheless that property (\ref{eq:K2}) is trivially 
fulfilled for  convolution kernels
$K_\theta (t, x)= K_\theta (t-x)$, 
which correspond to the standard kernel estimators.
\end{remark}
\par
The next example describes a collection of kernels corresponding to the
single--index model.
\paragraph*{\underline{Example.}}
{\em
%Consider the single--index model (i). 
Let $K:\rR^d\to \rR$, $\int K(t) dt=1$,
$E$ be an orthogonal matrix with the first vector--column equal to $e$.
Define for all $h\in \rR^d_+$
\[
K_h(t)= \Bigl[\prod_{i=1}^d h_i\Bigr]^{-1} K\left(\frac{t_1}{h_1}, \ldots,
\frac{t_d}{h_d}\right).
\]
Denote ${\cal H}=
\{h\in \rR^d_+: h=(h_1, h_{\max}, \ldots, h_{\max}), 
h_1\in [h_{\min}, h_{\max}]\}$, where the 
bandwidth range $[h_{\min}, h_{\max}]$
is supposed to be fixed. 
The collection of the kernels 
corresponding to the single--index model is
\[
{\cal K}=\Bigl\{K_\theta(t, x)= K_h[E^T(t-x)],\;\;\theta=(E, h)\in 
\Theta={\cal E}\times {\cal H}\subset \rR^d\Bigr\},
\]
where ${\cal E}$ is the set of all $d\times d$ orthogonal matrices.
\par
Clearly, $M({\cal K})=\|K\|_1$ so that K0 is fulfilled if 
$\|K\|_1<\infty$. Assumption K1 is trivially fulfilled because
$K_\theta(t, x)=K_{\theta}(t-x)$.
}
\par
For $\theta, \nu \in \Theta$ we define
\begin{equation}\label{eq:convoluted-kernel}
K_{\theta, \nu} (t, x) := \int K_\theta (t, y) K_\nu (y, x) dy,
\end{equation}
and let
\[
\hat{F}_{\theta, \nu} (x):=
\int K_{\theta, \nu}(t, x) Y(dt),\;\;\;x\in {\cal D}_0.
\]
Observe that $K_{\theta, \nu} = K_{\nu, \theta}$ in view of (\ref{eq:K2}),
so that indeed $\hat{F}_{\theta, \nu}\equiv \hat{F}_{\nu, \theta}$. 
This property is heavily exploited in the sequel, since 
the statistic $\hat{F}_{\theta, \nu}$ is an auxiliary estimator used in our construction. 
We have 
\begin{eqnarray}
\hat{F}_{\theta, \nu}(x) - F(x) &=& 
\int K_{\theta, \nu}(t, x) F(t) dt - F(x)  + \e\int K_{\theta, \nu}(t, x) W(dt)
\nonumber
\\
&=:& B_{\theta, \nu}(x) + \e Z_{\theta, \nu}(x).
\label{eq:Z-theta-nu1}
\end{eqnarray}
\par
The next simple result is a basic tool for construction of 
our selection procedure.
\begin{lemma}\label{lem:bias}
Let Assumption K0 hold;
then for any $F\in \rL_2(\cD) \cap \rL_p(\cD)$
\begin{equation}\label{eq:consequence}
\sup_{\nu\in \Theta}\|B_{\theta, \nu} - B_\nu\|_p \leq M({\cal K}) \|B_\theta\|_p\,,\;\;\;\forall \theta
\in \Theta.
\end{equation}
\end{lemma}
\pr By definition of $B_{\theta, \nu}$, $B_\nu$ and by the Fubini theorem
\begin{eqnarray*}
B_{\theta, \nu} (x)  - B_\nu (x) &=& 
\int K_{\theta, \nu}(t, x) F(t) dt  - \int K_{\nu}(t, x) F(t) dt 
\\
&=& \int K_\nu (y, x) 
\Bigl[\int K_{\theta}(t, y)  F(t) dt - F(y)\Bigr] dy  
\\
&=& \int K_\nu (y, x) B_\theta (y) dy.
\end{eqnarray*}
The statement of the lemma follows from the
general theorem about boundedness of integral operators
on $\rL_p$--spaces
[see, e.g., \citeasnoun[Theorem~6.18]{folland}] and (\ref{eq:M}).
\epr
\subsection{Selection rule}  
In order to present the basic idea underlying construction 
of the selection rule
we first discuss the noise--free version ($\e=0$)
of the estimation problem.
\paragraph*{\underline{Idea of construction (ideal case $\e=0$).}}
%Note that 
In this situation 
\[
{\cal F}({\cal K})=\Bigl\{\hat{F}_{\theta}(\cdot)=\int K_{\theta}(t, \cdot) F(t) dt, \;
\forall \theta\in \Theta\Bigr\}.
\]
so that
%Thus 
$\hat{F}_\theta$ can be viewed as a kernel--type 
approximation (smoother) of $F$. 
Note 
%also 
that the risk ${\cal R}_p[\hat{F}_\theta; F]=\|\hat{F}_\theta - F\|_p=\|B_{\theta}\|_p$ represents the quality
of approximation. 
%We denote by 
Let $\hat{F}_{\theta_*}$ be a smoother 
from ${\cal F}({\cal K})$
with the minimal approximation error, i.e.
%where
\[
\theta_*={\rm arg}
\inf_{\theta\in \Theta} {\cal R}_p[\hat{F}_{\theta}; F].
\]
Suppose that ${\cal K}$ satisfies Assumptions K0 and K1. Based on this collection 
we want to select a smoother, say $\hat{F}_{\hat{\theta}}\in 
{\cal F}({\cal K})$, that is
"as good as" $\hat{F}_{\theta_*}$, i.e., the smoother satisfying 
$\rL_p$--oracle inequality (\ref{oracle}). 
\par
To select 
$\hat{\theta}$ we suggest the following rule
\[
\hat{\theta} = {\rm arg} \inf_{\theta\in \Theta}\; 
\{\sup_{\nu\in \Theta} 
\|\hat{F}_{\theta, \nu} - \hat{F}_{\nu}\|_p \}.
\]
Let us compute the approximation error of the selected
smoother $\hat{F}_{\hat{\theta}}$.
By the triangle inequality  
\begin{eqnarray}
\|B_{\hat{\theta}}\|_p = \|\hat{F}_{\hat{\theta}} - F\|_p
&\leq& \|\hat{F}_{\hat{\theta}} - \hat{F}_{\hat{\theta}, \theta_*}\|_p +
\|\hat{F}_{\hat{\theta}, \theta_*} - \hat{F}_{\theta_*}\|_p + 
\|\hat{F}_{\theta_*}-F\|_p
\nonumber
\\
&=&  \|B_{\hat{\theta}} - B_{\hat{\theta}, \theta_*}\|_p +
\|B_{\hat{\theta}, \theta_*} - B_{\theta_*}\|_p+ \|B_{\theta_*}\|_p.
\label{eq:triangle}
\end{eqnarray}
In view of Assumption~K1 and (\ref{eq:consequence})
the first term on the right hand side of (\ref{eq:triangle})
does not exceed $M({\cal K})\|B_{\theta_*}\|_p$. 
To bound the second term we use
the definition of $\hat{\theta}$ and (\ref{eq:consequence}):
\begin{eqnarray*}
\|B_{\hat{\theta}, \theta_*} - B_{\theta_*}\|_p 
&\leq& 
\sup_{\nu\in \Theta}
\|B_{\hat{\theta}, \nu} - B_{\nu}\|_p 
\\
&\leq&
\sup_{\nu\in \Theta}
\|B_{{\theta_*}, \nu} - B_{\nu}\|_p \leq  M({\cal K})\|B_{\theta_*}\|_p.
\end{eqnarray*}
Combining these bounds we obtain from (\ref{eq:triangle}) that
\[
{\cal R}_p[\hat{F}_{\hat{\theta}}; F] \leq (2M({\cal K})+1) 
\|B_{\theta_*}\|_p = (2M({\cal K})+1)
\inf_{\theta\in\Theta} {\cal R}_p[\hat{F}_{\theta}; F].
\]
Therefore in the ideal situation $\e=0$, 
the $\rL_p$--oracle inequality (\ref{oracle})
holds with \mbox{${\cal L}=2M({\cal K})+1$}.
\paragraph*{\underline{Example (continuation).}}
{\em
We suppose additionally that there exists a positive integer 
$l$ such that
\[
\int t^{\blk} K(t) dt =0,\;\;\;|\blk| =1, \ldots, l, 
\]
where $\blk=(k_1, \ldots, k_d)$
is the 
multi--index, $k_i\geq 0$, $|\blk|=k_1+\cdots +k_d$, 
$t^{\blk}= t_1^{k_1}\cdots t_d^{k_d}$
for $t=(t_1,. \ldots, t_d)$.
Let $e$ is the true direction vector in the model (i). After rotation
described by  the matrix $E$ for any $h\in {\cal H}$ we have
\[
\|B_{\theta_*}\|_p \leq \Bigl\|
\int K(u)[f(\cdot+h_1u)-f(\cdot)] du\Bigr\|_p.
\]
If there exists $0<\alpha<l+1$, $L>0$ such that $f\in \rH_1(\alpha, L)$
then 
\begin{equation}\label{eq:bias-SI}
\|B_{\theta_*}\|_p\leq Lh_1^\alpha,\;\;\;\forall h_1\in [h_{\min}, h_{\max}].
\end{equation}
It is evident that when there is no noise in the model, the best choice of $h_1$
is $h_{\min}$.
}
\paragraph*{\underline{Idea of construction (real case $\e>0$).}}
When the noise is present, we use the same selection 
procedure 
with additional
control of  the noise contribution by its maximal value.
Similarly to the ideal case 
our selection rule is based on the statistics
$\{\sup_{\nu\in \Theta}\|\hat{F}_{\theta, \nu}-\hat{F}_\nu\|_p,
\;\theta\in \Theta\}$. Note that
\begin{eqnarray}
\|\hat{F}_{\theta, \nu}- \hat{F}_\nu\|_p  &\leq& 
\|B_{\theta, \nu} - B_\nu\|_p + \e\|Z_{\theta, \nu}-Z_\nu\|_p
\nonumber
\\
&\leq & 
\|B_{\theta, \nu} - B_\nu\|_p + 
\sup_x|\tilde{\sigma}_{\theta, \nu}(x)|\, 
\sup_{\theta, \nu} \|\tilde{Z}_{\theta, \nu}\|_p,
\label{eq:statistic}
\end{eqnarray}
where 
$Z_{\theta, \nu}(\cdot)$ and $Z_\nu (\cdot)$ are 
given 
in (\ref{eq:Z-theta-nu1}) and (\ref{eq:Z-theta}) respectively, and
\begin{eqnarray}
\sigma^2_{\theta, \nu}(x) &:=& \rE 
|Z_{\theta, \nu}(x)-Z_{\nu}(x)|^2=
\|K_{\theta, \nu}(\cdot, x) - K_{\nu}(\cdot, x)\|_2^2,
\;
x\in {\cal D}_0,
\label{eq:sigma-theta-nu}
\\
\tilde{\sigma}_{\theta, \nu}(x) &:=& \max\{\sigma_{\theta, \nu}(x)\,,\, 1\}
\nonumber
\\
\tilde{Z}_{\theta, \nu}(x)&:=&\tilde{\sigma}^{-1}_{\theta, \nu}(x)
[Z_{\theta, \nu}(x)-Z_{\nu}(x)].
\label{eq:Z-theta-nu}
\end{eqnarray}
\begin{remark}
In what follows we will be interested in large deviation probability
for the maximum of the  process $Z_{\theta, \nu}(x)-Z_{\nu}(x)$.
Typically the variance $\sigma_{\theta, \nu}(x)$ of this process 
tends to infinity as $\e\to 0$; therefore in the most interesting
examples $\tilde{\sigma}_{\theta, \nu}(x)=\sigma_{\theta, \nu}(x)$,
and $\tilde{Z}_{\theta, \nu}(x)$ has unit variance. However, for an abstract collection of the kernels, it can happen that $\sigma_{\theta, \nu}(x)$
is very small, for example, if $K_\theta$ approaches the delta--function. 
That is why we truncate the variance from below by $1$. 
%In this case the standard technique fails, and therefore we truncate 
%the variance from below by $1$. Then the problem is reduced
%to  the study of the non--standardized process,
%and the needed results are easily obtained.
\end{remark}
\par
In the ideal case we deduced from (\ref{eq:consequence}) 
that
\begin{equation}\label{eq:lower-estimator}
[M({\cal K})]^{-1} \sup_{\nu\in \Theta} \|\hat{F}_{\theta, \nu}- \hat{F}_\nu\|_p \leq \|B_{\theta}\|_p,\;\;\;\forall \theta\in \Theta,
\end{equation}
i.e., the left hand side can be considered as a lower estimator of
the bias. In the case of $\e>0$ 
we would like to guarantee the same property with high probability.
\par
This leads to the following control of the stochastic term.
Let $\delta\in (0,1)$, and  let $\kappa_p=\kappa_p({\cal K}, \delta)$ 
be the minimal positive real number
such that
\begin{equation}\label{eq:kappa}
\rP\Bigl\{\, 
\sup_{\theta\in \Theta}\; 
 \|\tilde{Z}_\theta(\cdot)\|_p
\geq \kappa_p\Bigr\}
\;+\;
\rP\Bigl\{\,
\sup_{(\theta, \nu) \in \Theta\times\Theta}\;\| 
\tilde{Z}_{\theta, \nu}(\cdot)\|_p\geq \kappa_p
\Bigr\} \leq \delta,
\end{equation}
where
similarly to (\ref{eq:sigma-theta-nu}) and (\ref{eq:Z-theta-nu})
we set
\begin{eqnarray*}
\tilde{Z}_\theta (x) &:=& \sigma^{-1}_{\theta}(x)Z_{\theta}(x),
\\
\sigma^2_\theta(x) &:=& \rE |Z_\theta(x)|^2 =\|K_\theta (\cdot, x)\|_2^2.
\end{eqnarray*}
The constant $\kappa_p$ controls deviation
of $\|\tilde{Z}_{\theta, \nu}\|_p$
as well as the deviation of standardized stochastic terms of all
estimators from the collection ${\cal F}({\cal K})$.  
We immediately obtain from (\ref{eq:statistic}), (\ref{eq:sigma-theta-nu})
and (\ref{eq:kappa}) that 
\begin{equation}\label{eq:b}
\hat{B}_\theta (p) := [M({\cal K})]^{-1}\sup_{\nu \in \Theta} 
 \big[\|\hat{F}_{\theta, \nu} - \hat{F}_\nu\|_p -
\e\kappa_p 
\sup_x \tilde{\sigma}_{\theta, \nu}(x)\big]
\leq \|B_\theta\|_p,\;\;\;\forall \theta\in \Theta,
\end{equation}
with probability larger than $1-\delta$. 
\par
Thus, similarly to (\ref{eq:lower-estimator}), 
$\hat{B}_\theta (p)$
is a lower estimator of the $\rL_p$--norm of the bias of the estimator
$\hat{F}_\theta$. This leads us to the following selection procedure.
\paragraph*{\underline{Selection rule.}} 
Define
\begin{equation}\label{eq:theta-hat}
\hat{\theta}= \hat{\theta}(\delta) := {\rm arg}\inf_{\theta\in \Theta}
\big\{\, \hat{B}_\theta(p) +
\kappa_p({\cal K}, \delta)\, \e\sup_x \sigma_\theta (x) \,\big\},
\end{equation}
and put finally
\[
\hat{F}(\delta)=\hat{F}_{\hat{\theta}}.
\]
\begin{remark}
The choice of $\hat{\theta}$ is very natural. Indeed, in view of
(\ref{eq:b}) for any $\theta\in\Theta$ with high probability
\begin{eqnarray*} 
\hat{B}_\theta(p) +
\kappa_p \e\sup_x \sigma_\theta (x) \leq \|{B}_\theta\|_p +
\kappa_p \e\sup_x \sigma_\theta (x).
\end{eqnarray*}
On the other hand, under rather general assumptions 
(see Section~\ref{subsec:L-infty})
\[
\|B_\theta\|_p+\e\kappa_p\sup_x 
\sigma_\theta (x)\leq C{\cal R}_p[\hat{F}_\theta; F],
\]
where $C$ is an absolute constant, independent of $F$ and $\e$.
Therefore with high probability
\[
\hat{B}_{\hat{\theta}}(p) +
\kappa_p \e\sup_x \sigma_{\hat{\theta}} (x) \leq 
C\inf_{\theta\in\Theta}
{\cal R}_p[\hat{F}_\theta; F].
\]
Thus in order to establish the $\rL_p$--norm oracle inequality it
suffices to majorate the risk of the estimator $\hat{F}_{\hat{\theta}}$
by $\hat{B}_{\hat{\theta}}(p) +
\kappa_p \e\sup_x \sigma_{\hat{\theta}} (x)$
and to choose $\delta=\delta(\e)$
tending to zero at an appropriate rate.
\end{remark}
\subsection{Basic result}
The next theorem establishes the basic result of this paper.
\begin{theorem}\label{th:main}
Let Assumptions K0 and K1 hold, and suppose that
\begin{itemize}
\item[(I)] $\hat{\theta}$ defined in (\ref{eq:theta-hat})
is measurable with respect to the observation $\{Y(t), t\in {\cal D}\}$,
and $\hat{\theta}$ belongs to $\Theta$;
\item[(II)] the events in (\ref{eq:kappa}) belong to the $\sigma$--algebra
generated by the observation $\{Y(t), t\in {\cal D}\}$.
\end{itemize}
Let $\delta\in (0,1)$, $\kappa_p$ be 
defined
in (\ref{eq:kappa}), and  $F$ be such that (I) and (II) hold. 
Then 
\begin{eqnarray}\label{eq:main-1}
\rE_F\|\hat{F}(\delta)-F\|_p  
\leq [3+2M({\cal K})] \inf_{\theta\in \Theta}
\big\{\,\|B_\theta\|_p + \kappa_p \e \big[\sup_x \sigma_\theta (x)\big]\big\} + r(\delta),
\end{eqnarray}
where
\[
r(\delta) := \|F\|_\infty [1+M({\cal K})] \delta + 
\sigma({\cal K}) \delta^{1/2} [\rE|\zeta|^2]^{1/2},
\]
$\sigma ({\cal K})$ is defined in (\ref{eq:sigma}), 
$\zeta:= \sup_{x, \theta} |\tilde{Z}_\theta (x)|$, and $\rE$ denotes
expectation with respect to the Wiener measure.
\end{theorem}
\begin{remark}\label{rem:T1}
In order 
to verify  
measurability of $\hat{\theta}$ and the condition (II)
we need to impose 
additional
assumptions on the collection of kernels ${\cal K}$. These assumptions
should guarantee   smoothness properties of the sample paths
of Gaussian processes 
$\{\tilde{Z}_\theta (x), (x, \theta)\in {\cal D}_0\times \Theta\}$ and
$\{\tilde{Z}_{\theta, \nu}(x), 
(x, \theta, \nu)\in \cD_0\times \Theta
\times \Theta\}$.
It is well--known [see, e.g., \citeasnoun{lifshits}] that such properties
for Gaussian processes can be described in terms of their covariance 
structures. In our particular case, the covariance structure is 
entirely determined  by the collection of kernels ${\cal K}$.
These fairly general conditions on ${\cal K}$ are given in 
Section~\ref{subsec:L-infty}.
\par
To ensure that $\hat{\theta}\in \Theta$ we need not only
smoothness conditions on the stochastic processes involved in the 
procedure description, but  also conditions on smoothness of $F$.
It is sufficient to suppose that $F$ belongs to some isotropic
H\"older ball, and this will be always assumed in the sequel.
This hypothesis
also guarantees that $F$ is uniformly bounded, 
which, in turn, implies boundedness of the remainder term $r(\delta)$. 
It is important to
note that neither procedure nor inequality (\ref{eq:main-1})
depend on parameters of this ball. 
\end{remark}
\begin{remark}
Our procedure 
% as well as 
and the basic oracle inequality
depend on the design parameter $\delta$.
% to be chosen.
The choice of this parameter is a delicate problem. 
%Indeed,
On the one hand, in order to reduce the remainder term we should choose
$\delta$ 
%should be 
as small as possible. On the other hand,
in view of the definition, $\kappa_p=\kappa_p(\delta)\to \infty$
as $\delta \to 0$.  
Note that we cannot minimize the right hand side 
of (\ref{eq:main-1}) with respect to $\delta$ 
because
%it would 
this leads to $\delta$ depending on unknown function
$F$. Fortunately, the same assumptions from Section~\ref{subsec:L-infty}
ensure that up to an absolute constant
\begin{equation}\label{eq:approx}
\inf_{\theta\in \Theta}
\big\{\,\|B_\theta\|_p + \kappa_p \e \big[\sup_x \sigma_\theta (x)\big]\big\}
\gtrsim  \e.
\end{equation}
The form of the remainder term $r(\delta)$ together with (\ref{eq:approx})
suggests 
% the choice of $\delta$
% dependent on $\e$, 
that $\delta$ should depend on $\e$,
for example,  
$\delta=\delta(\e)=\e^a$, $a>1$. Such a choice 
% of $\delta$ along with
under
assumptions from Section~\ref{subsec:L-infty} allows to show that
\begin{eqnarray}\label{eq:kappa-delta}
\kappa_p(\delta)=\kappa_p(\delta(\e))=\left\{\begin{array}{ll}
C(p), & p\in [1, \infty),\\
\sqrt{C(\infty)\ln (1/\e)}, & p=\infty,
                                             \end{array}
\right.
\end{eqnarray}
where $C(p)$, $p\in [1, \infty]$, are absolute constants, independent of
$\e$.
\end{remark} 
\par
Although the inequality (\ref{eq:main-1}) is not stated 
in the form of 
the $\rL_p$--norm oracle inequality, it can be helpful
(in view of (\ref{eq:kappa-delta})) 
 for deriving 
adaptive minimax results. To demonstrate 
%it 
this we return to the single--index model.
\paragraph{\underline{Example (continuation).}} 
{\em Remind that 
$\theta=(E, h)$ and note that
\begin{eqnarray*}
\sigma_\theta^2(x) = \sigma_{E, h}^2(x) = [h_1 h_{\max}^{d-1}]^{-2}
\int K_h^2[E^T(t-x)] dt
= [h_1 h_{\max}^{d-1}]^{-1}  \|K\|_2^2
\end{eqnarray*}
does not depend on  $E$ and $x$. 
% Let us 
Fix
$\delta=\e^a$ and 
%denote 
let $\hat{F}_\e$ be the estimator
$\hat{F}(\e^a)$ satisfying (\ref{eq:main-1}). 
Then (\ref{eq:main-1}) takes the form
\begin{eqnarray*}
\rE_F \|\hat{F}_\e - F\|_p &\leq& (3+2\|K\|_1) \inf_{E, h}
\big[ \|B_{E, h}\|_p+ \e\kappa_p(\e^a)  \sup_x \sigma_{E, h}(x)\big] +O(\e^a)
\\
&\leq&
(3+2\|K\|_1) \inf_h 
\big[\inf_E \|B_{E, h}\|_p + \e\kappa_p(\e^a)  [h_1 h_{\max}^{d-1}]^{-1/2}  \|K\|_2\big] +O(\e^a)
\\
&\leq &(3+2\|K\|_1)
\inf_{h_1}\big[ Lh_1^\alpha + \kappa_p(\e^a)  [h_1 h_{\max}^{d-1}]^{-1/2}  \|K\|_2\big] +O(\e^a).
\end{eqnarray*}
The last inequality follows from (\ref{eq:bias-SI}).
% To obtain the last inequality we used (\ref{eq:bias-SI}).
Taking into account (\ref{eq:kappa-delta}),
choosing $h_{\max}>0$ independent of $\e$, $h_{\min}=\e^2$,
and minimizing the last inequality 
with respect to $h_1\in [h_{\min}, h_{\max}]$ we obtain for all
$\alpha>0$, $L>0$
\[
\sup_{f\in \rH_1(\alpha, L)} \rE_F\|\hat{F}_\e - F\|_p
\leq C_p(L, h_{\max}, K)
\left\{\begin{array}{ll}
\e^{2\alpha/(2\alpha+1)}, & p\in [1, \infty)
\\*[2mm]
\big[\e\sqrt{\ln{(1/\e)}}\,\big]^{2\alpha/(2\alpha+1)}, & p=\infty.
       \end{array}
\right.
\]
It remains to 
note that $\hat{F}_\e$ does not depend on $(\alpha, L)$, and 
attains in view of 
the last inequality the minimax rate of convergence for all values
of
$(\alpha, L)$ simultaneously. It means that $\hat{F}_\e$ is optimally
adaptive on the scale of H\"older balls.
}
\subsection{Key oracle inequality}\label{subsec:L-infty}
In this section we discuss the choice of 
$\delta$ which leads to the {\em key oracle inequality}.
This inequality is suitable for deriving 
minimax and minimax adaptive results with minimal technicalities.
In particular, we  will use it for adaptive estimation in 
the additive multi--index model.
\par
In order to establish the {\em key oracle inequality} 
we need to impose
additional conditions on the collection of kernels ${\cal K}$. 
In particular, these conditions should guarantee the bounds 
(\ref{eq:kappa-delta})
for $\kappa_p(\delta(\e))$.
In the case $p=\infty$ such conditions are rather mild and standard;
they are related to 
deviation of supremum of Gaussian processes and therefore can be expressed
through  smoothness  of their covariance functions
\cite{lifshits}. As for the case $p<\infty$, we need to establish
bounds on large deviation probabilities of the 
$\rL_p$--norm of Gaussian  processes.
It requires additional assumptions on the collection of the kernels. Moreover,
such bounds cannot be directly obtained from the existing results.
We note nevertheless that  (\ref{eq:kappa-delta}) for the case $p<\infty$
can be  shown 
under fairly general assumptions, and this will be the subject of
a forthcoming paper.  From now on we restrict ourselves with the case
$p=\infty$.
\par
In the end of this section we discuss the connection between the 
{\em key
oracle inequality} and 
the {\em $\rL_\infty$--norm oracle inequality}
of type (\ref{oracle}).
\paragraph*{\underline{Assumptions.}} 
We suppose that the set $\Theta$ has the
following structure.
\begin{itemize}
\item[{\bf (A)}] $\Theta= \Theta_1 \times \Theta_2$ where
$\Theta_1=\{\theta^1, \ldots, \theta^N\}$ is a finite set, and
$\Theta_2\subset \rR^m$ is a compact subset of $\rR^m$ contained
in the Euclidean ball of radius $R$. Without loss of generality
we assume that $R\geq 1$.
\end{itemize}
\begin{remark}
Assumption~A allows to consider both discrete and continuous parameter
sets.
In particular, the case of empty $\Theta_2$ 
%can be viewed as the 
corresponds to 
selection from a finite set of estimators. This setup is often
considered
within the framework of the oracle approach.
In order to emphasize dependence of kernels $K_\theta$
on $\theta_1\in \Theta_1$ and $\theta_2\in \Theta_2$, we  sometimes write
$K_{(\theta_1, \theta_2)}$ instead of  $K_{\theta}$.
\end{remark}
\begin{itemize}
\item[{\bf (B)}] There exists $M_0$ such that $F\in \rH_d(M_0)$, where 
\[
\rH_d(M_0)=\Bigl\{g:\;\; g\in \bigcup_{\alpha>0, L>0} \rH_d(\alpha, L),\;
\|g\|_\infty\leq M_0\Bigr\}.
\]
\end{itemize}
\begin{remark}
Assumption~B is necessary 
% to justify 
for verification of the condition (I) of
Theorem~\ref{th:main}. It is also needed for deriving
the key oracle inequality from Theorem~\ref{th:main}
since it allows to bound uniformly the remainder term in (\ref{eq:main-1}). 
\par
We emphasize that our procedure does not depend on $M_0$. 
Finally note that
$\rH_d(M_0)$ is 
%very 
a huge set of functions(a bit smaller than the space of all bounded
continuous functions), i.e., Assumption~B is not restrictive at all.
\end{remark}

\begin{itemize}
\item[{\bf (K2)}] 
Denote $U:={\cal D}_0\times \Theta_2$.
There exist positive constants $\bar{L}$, 
and $\gamma\in (0,1]$ such that 
\begin{eqnarray*}
\hskip-1cm&&
\sup_{\theta_1\in \Theta_1}\;\sup_{u, u^\prime\in U}\; 
\frac{\|K_{(\theta_1, \theta_2)}(\cdot, x) - K_{(\theta_1, \theta_2^\prime)}
(\cdot, x^\prime)\|_2} 
{|u-u^\prime|^\gamma}\leq 
\bar{L},
\end{eqnarray*}
where $u=(x, \theta_2)$, and $u^\prime=(x^\prime, \theta_2^\prime)$. 
Without loss of generality we assume that $\bar{L}\geq 1$.
\end{itemize}
\begin{remark}
Assumption K2 ensures that sample
paths of the processes $\{\tilde{Z}_{\theta}(x), (x, \theta) \in \cD_0\times \Theta\}$ and 
$\{\tilde{Z}_{\theta, \nu}(x),  (x, \theta, \nu)\in \cD_0\times \Theta
\times \Theta\}$
belong with probability one 
to the isotropic  H\"older spaces $\rH_{m+d}(\tau)$ and
$\rH_{2m+d}(\tau)$ with regularity index $0<\tau<\gamma$
\cite[Section~15]{lifshits}. In particular, it is sufficient for 
fulfillment of
conditions (I) and (II) of Theorem~\ref{th:main}.
\end{remark}
\paragraph*{\underline{Choice of $\delta$.}}
Now we are ready to state the upper bound on the risk of our 
estimator (\ref{eq:theta-hat}) under Assumptions~A, B,~K0--K2.
Define
\[
C_{\cal K}:= M({\cal K}) \bar{L} R
\]
\begin{theorem}\label{th:main-2}
Let  Assumptions A, B,  K0--K2 hold, and assume that 
there exists $a>0$ such that
\begin{equation}\label{eq:choice-delta}
\delta_*:=
\min\Bigl\{ \frac{1}{N}, C_{\cal K}^{-(2m+d)/\gamma}, \e^2[\sigma({\cal K})]^{-2}\Bigr\} \geq \e^a.
\end{equation}
Let $\hat{F}_*=\hat{F}(\delta_*)$ be the estimator of Section~\ref{sec:scheme}
associated with the choice $\delta=\delta_*$. Then there exists a
constant $C_1\geq M_0$ depending on $d$, $m$ and $\gamma$ only such that 
\begin{eqnarray}
\rE_F \|\hat{F}_*-F\|_\infty
\leq [3+2M({\cal K})]
\inf_{\theta\in\Theta}\Bigl\{ \|B_{\theta}\|_\infty
+ C_1 \e\sqrt{\ln \e^{-1}} \sup_x \sigma_\theta(x)\Bigr\}.
\label{eq:main-2}
\end{eqnarray}
\end{theorem}
\begin{remark}
Typically in nonparametric setups $\bar{L}\sim \e^{-a_1}$,
$\sigma({\cal K})\sim \e^{-a_2}$ for some $a_1, a_2>0$.
If $N$ grows not faster than  $\e^{-a_3}$, then (\ref{eq:choice-delta})
holds.
\end{remark}
\paragraph*{\underline{$\rL_\infty$--norm oracle inequality.}}
Finally we show how the $\rL_\infty$--norm  
oracle inequality (\ref{oracle})
can be obtained from Theorem~\ref{th:main-2}.

\begin{theorem}\label{th:L-infty-oracle}
Assume that there exists a constant $C_2>0$ such that
\begin{equation}\label{eq:infty}
\inf_{\theta\in \Theta} \rE \|\tilde{Z}_\theta (\cdot)\|_\infty \geq
C_2\sqrt{\ln (1/\e)},
\end{equation}
and let $\hat{F}_*$ be the estimator from Theorem~\ref{th:main-2}.
Then
\begin{equation*}
 {\cal R}_\infty [\hat{F}_*;F] \leq {\cal L}\inf_{\theta\in\Theta}
 {\cal R}_\infty [\hat{F}_{\theta};F],
 \end{equation*}
where  ${\cal L}=[3+2M({\cal K})] \max\{1, C_1/C_2\}$.
\end{theorem}
\begin{remark}
The condition (\ref{eq:infty})
seems to be necessary
in order to have the constant 
${\cal L}$ independent of $\e$. 
In fact,  (\ref{eq:infty}) is an assumption on the collection of
kernels ${\cal K}$. To verify this condition one can use the Sudakov
lower bound on the expectation of the maximum of a Gaussian process
[see, e.g., \citeasnoun[Section~14]{lifshits}].
\end{remark}
\par
The proof of Theorem~\ref{th:L-infty-oracle} is an immediate
 consequence 
of Theorem~\ref{th:main-2}, (\ref{eq:infty}), and
the following auxiliary result that 
is interesting in its own right.
\begin{lemma}\label{lem:homotopy}
Let $\tilde{F}(\cdot)=\int S(t, \cdot) Y(dt)$ 
be a linear estimator of $F(\cdot)$. 
Denote by $B_S(\cdot)$ 
and $\e Z_S(\cdot)$
the bias and the stochastic part  
of $\tilde{F}(\cdot)-F(\cdot)$ respectively.
Then
for any  $F\in \rL_p({\cal D})\cap \rL_2({\cal D})$ and 
$p\in [1, \infty]$
\begin{eqnarray}\label{eq:risk-bound}
\frac{1}{4} \{ \|B_S\|_p + \e
\rE \|Z_S\|_p\}
\leq {\cal R}_p [\tilde{F}; F]
\leq  \|B_S\|_p + \e
\rE \|Z_S\|_p.
\end{eqnarray}
\end{lemma}
%

%%%%%%%%%%%%%%%%%%%%%%%
\section{Adaptive estimation in additive multi--index model}
\label{sec:multiindex}
In this section we apply the 
{\em key oracle inequality}  of Theorem~\ref{th:main-2}
to adaptive
estimation in the additive multi--index model.
\subsection{Problem formulation}
We impose that following structural
assumption on the function $F$ 
% observed 
in the model~(\ref{model}).
\par
Let ${\cal I}$ denote the set of all partitions of
$(1, \ldots, d)$, and for $\eta>0$ let 
\[
{\cal E}_\eta=\{ E=(e_1, \ldots, e_d): e_i \in \rS^{d-1},\;\;\;
|{\rm det}(E)|\geq \eta \}.
\]
For any $I\in {\cal I}$ and $E\in {\cal E}_\eta$ let
$E_1, \ldots, E_{|I|}$ be the corresponding partition of columns of $E$.
\begin{itemize}
\item[{\bf (F)}] 
{\em 
Let  
$I=(I_1, \ldots, I_{|I|})\in {\cal I}$, and $E\in {\cal E}_\eta$.
There exist
functions 
$f_i:\rR^{|I_i|}\to \rR$, $i=1, \ldots, |I|$ 
such that
\begin{equation*}
F(t) = \sum_{i=1}^{|I|} 
f_i( E^T_i t ).
\end{equation*}
}
\end{itemize}
\par
Assumption~F states  
that the unknown function $F$ can be represented as a sum 
of~$|I|$ unknown functions $f_i$, $i=1, \ldots, |I|$, where $f_i$ is
$|I_i|$--dimensional after an unknown linear
transformation. Note that partition $I$ is also unknown. The assumption
that $|{\rm det}(E)|\geq \eta$ is chosen for technical reasons; note 
%nevertheless 
that our estimation procedure does not require knowledge
of the value of this parameter. 
\par
Later on 
the functions $f_i$ 
will be 
are supposed to be smooth; in particular, we will assume
that all $f_i$'s belong to an isotropic H\"older ball
(see the next definition).
%
% The smoothness
%of each function, which will be also unknown, is understood as
%the belonging to an isotropic H\"older ball (see definition below).
\begin{definition}\label{def:holder}
A function $f:{\cal T} \to \rR$, ${\cal T}\subset \rR^s$, is said 
to belong to the 
H\"{o}lder ball $\rH_s(\beta, L)$ if $f$ has continuous partial
derivatives of all orders $\leq l$ satisfying the H\"{o}lder
condition with exponent $\alpha\in (0,1]$:
\begin{eqnarray*}
&&
\|D^{\blk} f\|_\infty \leq L,\;\;\; \forall |\blk|=0, \ldots, l;
\\*[2mm]
&&
\Bigl| f(z) - \sum_{j=0}^l \frac{1}{j!}
\sum_{|\blk|=j} D^{\blk} f(t) (z-t)^{\blk}\Bigr|
\leq L|z-t|^{\beta},\;\;\;
\forall z, t\in {\cal T},
\end{eqnarray*}
where $\beta=l+\alpha$, $\blk = (k_1, \ldots, k_s)$ is the 
multi--index, $k_i\geq 0$, $|\blk|=k_1+\cdots +k_s$, 
$t^{\blk}= t_1^{k_1}\cdots t_s^{k_s}$
for $t=(t_1,. \ldots, t_s)$, and 
$D^{\blk}=\partial^{|\blk|}/\partial t_1^{k_1}\cdots \partial t_s^{k_s}$.
\end{definition}
\par
The described  structure includes models (i)--(iv).
\begin{enumerate}
\item\  {\em [Single--index model.]} Let $F(t)=f(e^Tt)$ for some unknown
$f:\rR^d\to\rR$ and $e\in \rS^{d-1}$. In order to express the single--index 
model
in terms of assumption F, 
we set  $E=(e_1,\ldots, e_d)$ with
$e_1, e_2, \ldots, e_d$ being  an orthogonal basis of $\rR^d$ such that
$e_1=e$.  In this case we can set $I=(I_1, I_2)$ with $I_1=\{1\}$, 
$I_2=\{2, \ldots, d\}$ 
and $f_1=f$, $f_2\equiv 0$.
\item \ {\em [Additive model.]} Let $F(t)=\sum_{i=1}^d f_i(x_i)$
for unknown $f_i:\rR^d \to \rR$. Here 
$E$ is the $d\times d$ identity matrix, and $I=(I_1, \ldots, I_d)$, 
$I_i=\{i\}$.
\item \ {\em [Projection pursuit model.]}
Let $F(t)=\sum_{i=1}^d f_i(e_i^T t)$ for unknown $f_i:\rR^d\to \rR^1$
and unknown linearly independent 
direction vectors $e_1, \ldots, e_d\in \rS^{d-1}$.
Here 
$E=(e_1, \ldots, e_d)$, $I=(I_1, \ldots, I_d)$, $I_i=\{i\}$.
\item \ {\em [Multi--index model.]} Let 
$F(t)=f(e_1^Tt, \ldots, e_m^T t)$ 
for unknown direction vectors
$e_1, \ldots, e_m\in \rS^{d-1}$, and unknown function $f:\rR^m\to \rR^1$.
We define $E=(e_1, \ldots, e_d)$,
where $(e_{m+1}, \ldots, e_d)$ is the orthogonal basis
of the orthogonal complement to the subspace
${\rm span}\{e_1, \ldots, e_m\}$.
In this case we set $I=(I_1, I_2)$, $I_1=(1, \ldots, m)$, $I_2=(m+1, \ldots, d)$,
and $f_1=f$, $f_2\equiv 0$.

\end{enumerate}
\begin{definition}\label{def:F}
We say that function $F$ belongs to the class 
$\rF_{I, E}(\beta, L)$, $\beta>0$, $L>0$ if
\begin{itemize}
\item[(i)] Assumption~F is fulfilled with partition 
$I=(I_1, \ldots, I_{|I|}) \in {\cal I}$ and 
matrix $E\in {\cal E}_\eta$;
\item[(ii)] there exist positive real numbers $\beta_i$ and $L$
such that $f_i\in \rH_{|I_i|}(\beta_i, L)$, $i=1, \ldots, |I|$;
\item[(iii)] For all $i=1, \ldots, |I|$
\begin{equation}\label{eq:beta}
\beta=\frac{\beta_i}{|I_i|}.
\end{equation}
\end{itemize}
\end{definition}
\begin{remark}
The meaning of condition (iii) is that smoothness of functions $f_i$
is related to their dimensionality in such a way
that the effective smoothness of all functional components is the same. 
This condition 
does not restrict generality as  smoothness of a sum of functions
is determined by the worst smoothness of summands.
\end{remark}
\par
Let $\tilde{F}$ be an estimator of $F\in \rF_{I, E}(\beta, L)$; 
accuracy of $\tilde{F}$
is measured by the maximal risk
\[
{\cal R}_\infty 
[\tilde{F}; \rF_{I, E}(\beta, L)] := \sup_{F\in \rF_{I, E}(\beta, L)}
\rE_F \|\tilde{F} - F\|_\infty.
\]
\begin{proposition}[Minimax lower bound] 
\label{prop:lower}
Let $\varphi_\e(\beta)=[\e\sqrt{\ln(1/\e)}]^{2\beta/(2\beta+1)}$.
Then
\[
 \liminf_{\e\to 0} \inf_{\tilde{F}} \varphi_\e^{-1}(\beta)
\;{\cal R}_\infty 
[\tilde{F}; \rF_{I, E}(\beta, L)]>0,\;\;\;I\in {\cal I}, \;E\in {\cal E}_\eta,
\]
where $\inf$ is taken over all possible estimators $\tilde{F}$.
\end{proposition}
\begin{remark}
The appearance of the univariate rate 
$\varphi_\e(\beta)$
in the lower bound is not surprising since $2\beta/(2\beta+1)=2\beta_i/(2\beta_i+|I_i|)$, $i=1, \ldots, |I|$
 in view of (\ref{eq:beta}).
It is worth mentioning that 
$\varphi_\e (\beta)= \psi_{\varepsilon, |I_i|}(\beta_i)$
is the minimax rate of convergence in estimation of each
component $f_i$ [cf. (\ref{eq:psi})].
\end{remark}
The proof of Proposition~\ref{prop:lower} is absolutely standard 
and 
%can be 
is omitted. Obviously, the accuracy of estimation under
the additive multi--index model cannot 
be better than the accuracy of estimation of
one 
%each 
component provided that all other components are identically zero.
Since $E$ is fixed, the problem is reduced  to estimating 
$|I_i|$--variate function of smoothness $\beta_i$ in the model
(\ref{model}). In this case the lower bound 
is well--known and given by $\psi_{\varepsilon, |I_i|}(\beta_i)$. It remains to note that $\psi_{\varepsilon, |I_i|}(\beta_i)$ does not depend on $i$ and coincides
with $\varphi_\e(\beta)$ in view of (\ref{eq:beta}).
\par
Below we propose an estimator that attains the rate
$\varphi_\e(\beta)$ 
simultaneously over $\rF_{I, E}(\beta, L)$, $I\in {\cal I}$, 
$E\in {\cal E}_\eta$, $0<\beta\leq \beta_{\max}<\infty$, $L>0$, i.e.,
the optimally adaptive estimator.
\subsection{Kernel construction}
To construct a family of kernel estimators let us
consider the idealized situation  when both the partition 
$I=(I_1, \ldots, I_{|I|})\in {\cal I}$ 
and $E\in {\cal E}_\eta$ are known.
\par
\begin{itemize}
\item[{\bf (G)}]
{\em
 Let $g:[-1/2, 1/2] \to \rR$ be a univariate
kernel satisfying the following conditions
\begin{itemize}
\item[(i)]
$
\int g(x) dx =1,\;\;\;\int g(x) x^k dx=0,\;\;k=1, \ldots, \ell$;
\item[(ii)]  $g\in \rC^1$.
\end{itemize}
}
\end{itemize}
\par
Fix  a bandwidth $h=(h_1, \ldots, h_d)$, 
$h_{\min} \leq h_i\leq h_{\max}$
and put
\begin{eqnarray*}
G_0(t) &=& \prod_{i=1}^d g(t_i)
\\
G_{i, h}(t) &=& 
\prod_{j \in I_i} \frac{1}{h_j} g\Bigl(\frac{t_j}{h_j}\Bigr)
\prod_{j\not\in I_i} g(t_j),\;\;\
i=1, \ldots, |I|.
\end{eqnarray*}
Now we define the kernel associated with partition $I$,
matrix $E$, 
and bandwidth $h$. 
Fix  $\theta=(I,E, h)\in \Theta=
{\cal I}\times {\cal E}_\eta\times [h_{\min}, h_{\max}]^d$, and let
\begin{equation}\label{eq:kernel}
K_{\theta}(t) = |{\rm det}(E)|\sum_{i=1}^{|I|} G_{i,h}(E^Tt) - (|I|-1)
|{\rm det}(E)| G_0(E^Tt).
\end{equation}
\subsection{Properties of the kernel}
First we  
state
%formulate 
%the 
evident properties
of the kernel $K_\theta$.
\begin{lemma}\label{lem:var}
For any $\theta\in\Theta$ 
\begin{eqnarray}
&&
\int K_{\theta}(t) dt = 1
%\label{eq:int-K-1}
\nonumber
\\
&&
\|K_{\theta}\|_1 \leq  
(2|I|-1)\|g\|_1^{d}.
%\label{eq:K-norm-1}
\nonumber
\\
&&
\|K_{\theta}\|_2 \leq  
 |{\rm det(E)}|^{1/2} \|g\|_2^d
\Bigl( \sum_{i=1}^{|I|} 
\prod_{j\in I_i} h_j^{-1/2} + |I|-1\Bigr).
\label{eq:K-norm-2}
\end{eqnarray}
\end{lemma}
The proof follows straightforwardly from (\ref{eq:kernel}).
\par
Next lemma establishes approximation properties of $K_\theta$.
Put for any $x\in \cD_0$
\begin{eqnarray*}
B_\theta (x) &=& \int K_{\theta}(t-x) F(t)dt -F(x).
\end{eqnarray*}
Clearly, $B_\theta(\cdot)$ is
the 
bias of the estimator 
associated with kernel $K_\theta$.
\begin{lemma}\label{lem:approx}
Let $F\in \rF_{I, E}(\beta, L)$, and let Assumption~G hold with
$\ell=\max_i \lfloor \beta_i\rfloor$. Then
\begin{equation}\label{eq:bias}
\|B_{\theta}\|_\infty \leq L\sum_{i=1}^{|I|} 
\|g\|_1^{|I_i|} \sum_{j\in I_i} h_j^{\beta_i}.
\end{equation}
\end{lemma}
\begin{remark}
Lemmas~\ref{lem:var} and~\ref{lem:approx} allow to 
%explain
derive an upper bound on
the accuracy of estimation on the class $\rF_{I, E}(\beta, L)$
for given $I$ and $E$. Indeed, the typical balance equation
for the bandwidth selection takes the form
%for $\rL_\infty$--losses 
%
%permitting to choose the bandwidth in the
%optimal way is
\[
\e\sqrt{\ln (1/\e)}\|K_\theta\|_2= \|B_\theta\|_\infty.
\]
Therefore using the upper bounds in (\ref{eq:bias})
and (\ref{eq:K-norm-2}) we arrive to the optimal choice of bandwidth
given by $h=h^*=(h_1^*, \ldots, h^*_d)$,
\begin{equation}\label{eq:ideal-band-choice}
h^*_j= \Bigl(\frac{\e}{L}\sqrt{\ln (1/\e)}\Bigr)^{2/(2\beta_i+{|I_i|})}
\Bigl(\frac{\|g\|_2}{\|g\|_1}\Bigr)^{2d/(2\beta_i+{|I_i|})},\;\;\;j\in I_i,\;\;
i=1, \ldots, |I|.
\end{equation}
If $\hat{F}_\theta (x) = \int K_{\theta}(t-x) Y(dt)$
is a kernel estimator with $\theta=(I, E, h_*)$ then
we have the following 
upper bound on its $\rL_\infty$--risk:
\begin{equation}\label{eq:upper}
{\cal R}_\infty [\hat{F}_\theta; \rF_{I, E}(\beta, L)] \leq 
C L^{1/(2\beta+1)}\varphi_\e(\beta),
\end{equation}
where $C$ is an absolute constant. Thus, in view of Proposition~\ref{prop:lower},
$\varphi_\e(\beta)$
is the minimax rate of convergence on the class $\rF_{I, E}(\beta, L)$.
We stress that 
%the 
construction of minimax estimator $\hat{F}_\theta$
requires
%the 
knowledge of all parameters of the functional class: $I$, $E$, $\beta$
and~$L$.
\end{remark}

\subsection{Optimally adaptive estimator}
%Choose 
Let $h_{\min}=\e^2$ and $h_{\max}=\e^{2/[(2\beta_{\max}+1)d]}$ 
for some 
$\beta_{\rm max}>0$.
Consider the collection of kernels 
${\cal K}=\{K_\theta (\cdot),\,\theta =(I, E, h)\in \Theta\}$
where  $K_\theta(\cdot)$ is defined in (\ref{eq:kernel}). 
The corresponding collection of estimators is given by
\[
{\cal F}({\cal K})=\Bigl\{
\hat{F}_\theta (x) = \int K_{\theta}(t-x) Y(dt),\;\;\;\theta\in \Theta
\Bigr\}.
\]
Based on the collection ${\cal F}({\cal K})$ we define the estimator
$\hat{F}_*$ following the selection rule 
(\ref{eq:theta-hat}) with the choice of $\delta=\e^a$ 
where $a=24d^3+12d^2$.
\begin{theorem}\label{th:multi}
Suppose that Assumption G holds with $\ell= \lfloor d\beta_{\max} \rfloor$.
Then for any $I\in {\cal I}$, $E\in {\cal E}_\eta$, 
$0<\beta\leq \beta_{\max}$, and $L>0$ 
\begin{equation*}
\limsup_{\e\to 0}
\varphi_\e^{-1}(\beta)\,{\cal R}_\infty [\hat{F}_*; \rF_{I, E}(\beta, L)] \leq 
C L^{1/(2\beta+1)},
\end{equation*}
where $C$ depends on $d$, $\beta_{\max}$,  and 
the kernel $g$ only.
\end{theorem}
\par
Combining the results of Theorem~\ref{th:multi} and 
Proposition~\ref{prop:lower} we obtain that the estimator $\hat{F}_*$
is optimally adaptive on the scale of functional classes
$\rF_{I, E}(\beta, L)$.
Thus this estimator adjusts automatically to unknown structure
as well as to unknown smoothness.
\par
We note that
traditionally any structural assumption is understood as the existence
of the structure. Mathematically in our case it means that the underlying function
belongs to the union of classes $\rF_{I, E}(\beta, L)$ with respect to $I\in {\cal I}$
and $E\in {\cal E}_\eta$, i.e.,
\[
F\in \rF(\beta, L)=\bigcup_{I\in {\cal I}, E\in {\cal E}_\eta}
\rF_{I, E}(\beta, L).
\]
Next theorem shows that our estimation procedure
is optimally adaptive
on the scale of functional classes $\rF(\beta, L)$, $0<\beta\leq \beta_{\max}$,
$L>0$.
\begin{theorem}\label{th:multi-1}
Suppose that Assumption G holds with $\ell=\lfloor d\beta_{\max}\rfloor$.
Then for any  
$0<\beta\leq \beta_{\max}$, and $L>0$ 
\begin{equation*}
\limsup_{\e\to 0}
\varphi_\e^{-1}(\beta)\,{\cal R}_\infty [\hat{F}_*; \rF(\beta, L)] \leq 
C L^{1/(2\beta+1)},
\end{equation*}
where $C$ depends on $d$, $\beta_{\max}$,  and 
the kernel $g$ only.
\end{theorem}
Theorem~\ref{th:multi} follows immediately from Theorem~\ref{th:multi-1}.
Proposition~\ref{prop:lower} together with Theorem~\ref{th:multi-1}
shows that 
in terms of rates of convergence 
there is no price to pay for adaptation with respect to unknown
structure.

\section{Proofs of Theorems~\ref{th:main},~\ref{th:main-2} and~\ref{th:multi-1}}
\label{sec:proofs}
\paragraph*{Proof of Theorem~\ref{th:main}.}
Define the random event
\begin{equation*}
A= A_1 \cap A_2:=\Bigl\{\omega: 
\sup_{\theta \in \Theta} \|\tilde{Z}_\theta\|_p \leq \kappa_p
\Bigr\}
\cap 
\Bigl\{\omega: \sup_{(\theta, \nu)\in \Theta\times\Theta } 
\|\tilde{Z}_{\theta, \nu}\|_p\leq
\kappa_p\Bigr\}. 
\end{equation*}
\par
1$^0$.
First, we observe that
\begin{equation}\label{eq:b-bound}
\hat{B}_\theta(p)\, 1(A) \leq \|B_\theta\|_p,\;\;\;\forall \theta\in \Theta.
\end{equation}
Indeed, in view of Lemma~\ref{lem:bias} 
on the set $A$
\begin{eqnarray*}
\|B_{\theta}\|_p &\geq & \sup_{\nu\in \Theta}\frac{1}{\|K_\nu\|_{1, \infty}}
\Bigl\|\int K_\nu(t, x) B_\theta(t) dt \Bigr\|_p
\\
&\geq &  [M({\cal K})]^{-1}
\sup_{\nu\in\Theta} 
\Bigl(
\|\hat{F}_{\theta, \nu} - \hat{F}_\nu\|_p
- \e \|Z_{\theta, \nu} - Z_{\nu}\|_p \Bigr)
\\
&\geq & 
[M({\cal K})]^{-1}
\sup_{\nu\in\Theta} 
\big[
\|\hat{F}_{\theta, \nu} - \hat{F}_\nu\|_p
- \kappa_p\e \sup_x \tilde{\sigma}_{\theta, \nu}(x)\big]
= \hat{B}_\theta(p),
\end{eqnarray*}
where we have also used definition of $A$ and the fact that
\[
\hat{F}_{\theta, \nu}(x)-\hat{F}_\nu(x)= \int K_\nu (t, x) B_\theta(t) dt +
\e [Z_{\theta, \nu}(x)-Z_\nu(x)].
\]
\par
$2^0$. Second, we note that for any $\theta, \nu\in \Theta$
\begin{eqnarray}
\sup_x \sigma_{\theta, \nu}(x) &=&
\|K_{\theta, \nu} - K_\nu\|_{2, \infty} 
\leq  \|K_{\theta, \nu}\|_{2, \infty} +\|K_\nu\|_{2, \infty}
\nonumber
\\
&\leq& \|K_\theta\|_{1, \infty} \|K_\nu\|_{2, \infty} + \|K_{\nu}\|_{2, \infty}
\leq [1+M({\cal K})]\, \|K_\nu\|_{2, \infty}
\nonumber
\\
&=&  [1+M({\cal K})] \sup_x \sigma_\nu (x).
\nonumber
\end{eqnarray}
Here we have used the inequality $\|K_{\theta, \nu}\|_{2, \infty}\leq 
\|K_\theta\|_{1, \infty} \|K_\nu\|_{2, \infty}$ which follows from the 
Minkowski integral inequality.
\par
The Cauchy--Schwarz inequality and (\ref{eq:K1}) yield
$\sigma_\nu(x)\geq (\mes\{\cD\})^{-1/2}$ for all $x$ and $\nu$. 
This 
implies without loss of generality that for any 
$\theta, \nu\in \Theta$
\begin{equation}
\label{eq:ineq}
\sup_x \tilde{\sigma}_{\theta, \nu}(x) \leq [1+M({\cal K})] 
\sup_x \sigma_\nu (x).
\end{equation}
\par
$3^0$.
Now define
\[
\theta_* := {\rm arg} \inf_{\theta\in \Theta} \{\,\|B_\theta\|_p + \kappa_p \e
\sup_x \sigma_\theta (x)\,\},
\]
and let $\hat{F}_* = \hat{F}_{\theta_*}$.
We write
\begin{eqnarray}
\|\hat{F}-F\|_p 1(A) &\leq& \|\hat{F}_{\theta_*} - F\|_p\, 1(A)
+
\|\hat{F}_{\theta_*} - \hat{F}_{\hat{\theta}, \theta_*}\|_p\, 1(A)
\nonumber
\\
&&\hspace{30mm}
+\;\;
\|\hat{F}_{\hat{\theta}} - \hat{F}_{\hat{\theta}, \theta_*}\|_p\, 1(A),
\label{eq:error}
\end{eqnarray}
and note that
\begin{equation}\label{eq:first-term}
\|\hat{F}_{\theta_*} - F\|_p\, 1(A) \leq \|B_{\theta_*}\|_p  + \kappa_p\e
\sup_x\sigma_{\theta_*}(x) = 
\inf_{\theta\in\Theta} \{\, \|B_\theta\|_p +
\kappa_p\e\sup_x\sigma_\theta (x)\,\}.
\end{equation}
Furthermore, 
\begin{eqnarray}
\|\hat{F}_{\theta_*} - \hat{F}_{\hat{\theta}, \theta_*}\|_p\, 1(A)
&\leq & M({\cal K}) \hat{B}_{\hat{\theta}}(p)\, 1(A) + 
\kappa_p\e\sup_x \tilde{\sigma}_{\hat{\theta}, \theta_*}(x)
\nonumber
\\
&\leq& M({\cal K}) \hat{B}_{\hat{\theta}}(p)\, 1(A)
+  [1+M({\cal K})] \kappa_p\e \sup_x \sigma_{\theta_*}(x), 
\nonumber
\end{eqnarray}
where the first inequality follows from definition of $\hat{B}_\theta(p)$; 
the second inequality is a consequence of (\ref{eq:M}) and (\ref{eq:ineq}).
Similarly, 
\begin{eqnarray}
\|\hat{F}_{\hat{\theta}} - \hat{F}_{\hat{\theta}, \theta_*}\|_p 1(A) 
& \leq &
M({\cal K}) \hat{B}_{\theta_*}(p)\, 1(A)  + 
\kappa_p \e \sup_x \tilde{\sigma}_{\theta_*, \hat{\theta}}(x)
\nonumber
\\
&\leq&
M({\cal K}) \hat{B}_{\theta_*}(p)\, 1(A) + [1+M({\cal K})] \kappa_p \e 
\sup_x \sigma_{\hat{\theta}}(x).
%\label{eq:third-term}
\nonumber
\end{eqnarray}
Now using (\ref{eq:theta-hat}) and (\ref{eq:b-bound}) we obtain 
\[
\begin{array}{l}
[\|\hat{F}_{\theta_*} - \hat{F}_{\hat{\theta}, \theta_*}\|_p 
+
\|\hat{F}_{\hat{\theta}} - \hat{F}_{\hat{\theta}, \theta_*}\|_p] 1(A) 
\\*[2mm]
\hspace{10mm}\leq  
[1+M({\cal K})]
\big\{[\hat{B}_{\hat{\theta}}(p) + \hat{B}_{\theta_*}(p)]1(A) + 
\kappa_p \e\sup_x \sigma_{\hat{\theta}}(x) 
+ \kappa_p \e \sup_x \sigma_{\theta_*}(x)\big\}
\\*[2mm]
\hspace{10mm}
\leq  2[1+M({\cal K})] \big\{\|B_{\theta_*}\|_p + \kappa_p\e 
\sup_x \sigma_{\theta_*}(x)\big\}.
\end{array}
\]
Then (\ref{eq:error}) and (\ref{eq:first-term}) lead to
\begin{equation}\label{eq:good-set}
\|\hat{F}-F\|_p 1(A) \leq [3+2M({\cal K})]
 \inf_{\theta \in \Theta}
\{\,\|B_\theta\|_p + \kappa_p \e \sup_x \sigma_\theta (x)\,\}.
\end{equation}
\par
4$^0$. In order to complete the proof 
it suffices to bound $\|\hat{F}-F\|_p 1(A^c)$.
Note that by our choice of $\kappa_p$ (see (\ref{eq:kappa})),
$\rP(A^c)\leq \delta$.
Moreover
\begin{eqnarray*}
\|\hat{F}-F\|_p 1(A^c) &\leq & 
( \sup_{\theta\in \Theta}\|B_\theta\|_p  +
\sup_{\theta\in \Theta} \|Z_\theta(\cdot)\|_p )\, 1(A^c)
\\
&\leq & \|F\|_\infty [1+M({\cal K})]
1(A^c) + \sigma({\cal K}) \zeta\, 1(A^c),
\end{eqnarray*}
where $\sigma({\cal K})$ is defined in (\ref{eq:sigma}), and
$\zeta := 
\sup_{x, \theta} |\tilde{Z}_\theta (x)|$. Therefore
\begin{eqnarray*}
\rE \|\hat{F}-F\|_p 1(A^c)
 &\leq&  \|F\|_\infty [1+M({\cal K})] \rP (A^c) + \overline{\sigma} 
[\rE\zeta^2]^{1/2} \rP^{1/2}(A^c)
\\*[2mm]
&\leq &
\|F\|_\infty [1+M({\cal K})] 
\delta +  \sqrt{\delta}\,\bar{\sigma} [\rE|\zeta|^2]^{1/2}
\end{eqnarray*}
where 
we have used 
(\ref{eq:kappa}). 
Combining this inequality with (\ref{eq:good-set}) we
complete the proof.
\epr
\paragraph{Proof of Theorem~\ref{th:main-2}.}	
\par
1$^0$. First we show that
Assumptions A, B, and K2 imply conditions (I) and (II)
of Theorem~\ref{th:main}.
\par
Indeed, Assumption K2 ensures that sample
paths of the processes $\{\tilde{Z}_{\theta}(x), (x, \theta) \in \cD_0\times \Theta\}$ and 
$\{\tilde{Z}_{\theta, \nu}(x),  (x, \theta, \nu)\in \cD_0\times \Theta
\times \Theta\}$
belong with probability one 
to the isotropic  H\"older spaces $\rH_{m+d}(\tau)$ and
$\rH_{2m+d}(\tau)$ with regularity index $0<\tau<\gamma$
\cite[Section~15]{lifshits}. 
Thus the condition~(II) is fulfilled.
\par
Moreover,
together with Assumption~B this implies that
for any $F\in \rH_d(M_0)$
sample paths of the process
$\hat{F}_{\theta,\nu}(x)-\hat{F}_\nu (x)$ belong with probability 
one 
to the isotropic  H\"older space 
$\rH_{2m+d}(\tau^\prime)$ 
on $\cD_0\times\Theta\times\Theta$
with some regularity index $0<\tau^\prime<\gamma$. 
This, in turn, shows that  for any $F\in \rH_d(M_0)$ 
sample paths of the process
\[
\sup_{\nu\in \Theta}\|\hat{F}_{\theta,\nu}-\hat{F}_\nu\|_p
\]
belong to $\rH_{m}(\tau^\prime)$ on $\Theta$. 
Then condition~(I) holds
% This together with
in view of Assumption~A 
%yields condition (I) in view of
and \citeasnoun{jennrich}.
\par
2$^0$.
It follows from Lemma~\ref{lem:exp-inequalities} in Appendix that for any 
$\kappa\geq 1+\sqrt{(2m+d)/\gamma}$
\begin{eqnarray*}
&& \rP\Bigl\{ \sup_{\theta\in\Theta} \|\tilde{Z}_\theta (\cdot) \|_\infty 
\geq \kappa\Bigr\} 
+
\rP\Bigl\{ \sup_{(\theta, \nu)\in \Theta\times\Theta} \|\tilde{Z}_{\theta, \nu}(\cdot)\|_\infty 
\geq \kappa\Bigr\}
\\
&& \hspace{45mm} \leq N^2[c_1 M({\cal K})\bar{L}R \kappa]^{(2m+d)/\gamma} \exp\{-\kappa^2/2\},
\end{eqnarray*}
where $c_1$ is an absolute constant.
By definition of $\kappa$
we obtain that
\begin{equation}\label{eq:kappa-in}
\exp\{\kappa^2/2\}
\leq N^2[c_1 M({\cal K})\bar{L} R\kappa ]^{(2m+d)/\gamma}\delta_*^{-1}
\end{equation}
which, in turn, implies
\begin{eqnarray}
\kappa &\leq& \Bigl[2\ln \delta_*^{-1} + 4\ln N + \frac{2(2m+d)}{\gamma} \ln 
C_{{\cal K}} + \frac{2m+d}{\gamma} (\ln \kappa^2 + c_2)\Bigr]^{1/2}
\nonumber
\\
&\leq& \sqrt{c_3\ln \e^{-1}} =: \bar{\kappa},
\label{eq:kappa-upper}
\end{eqnarray}
where $c_3$ depends on $(2m+d)/\gamma$ only; here we have
used (\ref{eq:choice-delta}). 
\par
Now we bound the remainder term in (\ref{eq:main-1}).  
It follows from Lemma~\ref{lem:exp-inequalities} that
for any $\lambda\geq 1+ \sqrt{(d+m)/\gamma}$ one has
\begin{eqnarray*}
\rE|\zeta|^2 &=& \int_0^\infty 2t \rP(\zeta > t) dt
\leq
2\lambda + 
2 \int_{\lambda}^\infty t N[c_4\bar{L}R t]^{(d+m)/\gamma} e^{-t^2/2}
dt
\\
&\leq& 2\lambda + 2N [c_4\bar{L}R]^{(d+m)/\gamma} e^{-\lambda^2/4}
\int_0^\infty t^{1+(d+m)/\gamma} e^{-t^2/4} dt.
\end{eqnarray*}
If we choose $\lambda=\sqrt{2}\bar{\kappa}$ and apply
(\ref{eq:kappa-in}), we get
\[
\rE |\zeta|^2 \leq 2\sqrt{2}\bar{\kappa} + c_5 N^{-1}\delta_*
\leq c_6 \ln \delta_*^{-1}.
\]
Using  (\ref{eq:choice-delta}) and the fact that $\sigma({\cal K})\geq c_7$
 we finally obtain 
$r(\delta_*)\leq M_0 [1+M({\cal K})]\e+ c_8\e\sqrt{\ln\e^{-1}}$
which yields~(\ref{eq:main-2}).
\epr

\paragraph*{Proof of Theorem~\ref{th:multi-1}.}
1$^0$. In order to apply the result of Theorem~\ref{th:main-2} we have to verify
Assumption~K2 for the collection of kernels defined in (\ref{eq:kernel}).
Recall that $\theta=(I, E, h)$, and in notation of 
Assumptions~A and~K2, $\theta=(\theta_1, \theta_2)$, where 
$\theta_1=I\in \Theta_1={\cal I}$, and $\theta_2=(E, h)\in \Theta_2={\cal E}_\eta\times [h_{\min}, h_{\max}]^d$.
\par 
We deduce from (\ref{eq:kernel}) and Assumption~G(ii) that 
$K_\theta (t)$ is continuously differentiable in $\theta_2$ and $t$, and
\[
\sup_{\theta_2\in \Theta_2}\; \sup_{t\in \cD} 
|\nabla_{\theta_2, t}\; K_\theta (t)| 
\leq \tilde{L}h_{\min}^{-3d},
\]
where $\tilde{L}$ is an absolute constant depending only on $d$ and 
$\|g\|_\infty$. Taking into account that $h_{\min}=\e^2$ we arrive to
Assumption~K2 with 
\begin{equation}\label{eq:bar-L}
\bar{L}=\tilde{L}\e^{-6d},\;\;\;\hbox{and}\;\;\; \gamma=1/2.
\end{equation}
\par
2$^0$. In view of (\ref{eq:bar-L}), assumption (\ref{eq:choice-delta})
is verified.
\par
3$^0$. Fix $\beta$ and $L$ and assume that
$F\in \rF(\beta, L)$. By definition of the class  
$F\in \rF(\beta, L)$ there exist $I_*\in {\cal I}$ and $E_*\in {\cal E}_\eta$
such that $F\in \rF_{I_*, E_*}(\beta, L)$. Let $h_*$ be given
by (\ref{eq:ideal-band-choice}). Then from (\ref{eq:main-2})
and (\ref{eq:upper})
\begin{eqnarray*}
&&\rE_F \|\hat{F}_*-F\|_\infty
\\
&&\leq  [3+2M({\cal K})]
\inf_{(I, E, h)\in\Theta}\Bigl\{ \|B_{I, E, h}\|_\infty
+ C_1 \e\sqrt{\ln \e^{-1}} \sup_x \sigma_{I, E, h}(x)\Bigr\}
\\
&&\leq
[3+2M({\cal K})]
\Bigl\{ \|B_{I_*, E_*, h_*}\|_\infty
+ C_1 \e\sqrt{\ln \e^{-1}} \sup_x \sigma_{I_*, E_*, h_*}(x)\Bigr\}
\\
&&
\leq 2[3+2M({\cal K})] (C_1\vee 1)
C L^{1/(2\beta+1)}\varphi_\e(\beta),
\end{eqnarray*}
where $C$ is the constant appearing in (\ref{eq:upper}).
\epr
\section*{Appendix}
\paragraph*{Proof of Lemma~\ref{lem:homotopy}.}
Only the left hand side inequality should be
proved.
First we note that
\[
\|f\|_p = \sup\Bigl\{ |\int \phi f |: \;\|\phi\|_q=1\Bigr\} 
\]
\cite[p.~188]{folland}.
Thus we have for $p<\infty$
\begin{eqnarray*}
\rE_F \|\tilde{F}-F\|_p
&=& \rE_F \| B_S + \e Z_S\|_p 
\\
&=& \rE_F  \sup_{g: \|g\|_q\leq 1} \int [B_S(x) +\e Z_S(x)] g(x) dx
\\
&\geq & \rE_F \int [B_S(x)+ \e Z_S(x)] g_*(x) dx,
\end{eqnarray*}
where 
$g_*(x) = \|B_S\|_p^{-p/q} |B_S(x)|^{p-1} \sign\{B_S(x)\}$. 
Therefore  
\begin{eqnarray}\label{eq:B-1}
\rE_F \| B_S + \e Z_S\|_p \geq \int B_S(x) g_*(x) dx +
\rE \int Z_S(x) g_*(x) dx = \|B_S\|_p.
\end{eqnarray}
On the other hand, by the triangle inequality
$\rE_F \| B_S + \e Z_S\|_p \geq \e \rE \|Z_S\|_p - \|B_S\|_p$. 
Combining the two last inequalities we obtain 
$\rE_F \| B_S + \e Z_S\|_p \geq \frac{1}{2} \e \rE \|Z_S\|_p$ which along with
(\ref{eq:B-1}) yields (\ref{eq:risk-bound}).
\par
If $p=\infty$ then 
for any $x_0\in {\cal D}_0$ one has 
$
\rE\|B_{\theta}+\e Z_\theta\|_\infty \geq \pm \rE [B_{\theta}(x_0)+ \e Z_{\theta}(x_0)]
= \pm B_\theta (x_0)$, and therefore 
$\rE \|B_{\theta}+ \e Z_\theta\|_\infty \geq \|B_{\theta}\|_\infty$.
\epr
\paragraph*{Proof of Lemma~\ref{lem:approx}.}
We will use the following notation: 
for any vector $t\in \rR^d$, and partition $I=(I_1, \ldots, I_{|I|})$
we will write
$t_{(i)}=(t_j, j\in I_i)$.
Throughout the proof without loss of generality we assume that
$E$ is the $d\times d$ identity matrix.
\par
Using the fact that $F(t)=\sum_{i=1}^{|I|} f_i(E_i^Tt)$ 
we have
\begin{eqnarray*}
\int K_{\theta}(t-x) F(t) dt = 
\sum_{i=1}^{|I|} \sum_{j=1}^{|I|} \int G_{j,h}(t-x) f_i(t_{(i)}) dt
- (|I|-1) \sum_{i=1}^{|I|} \int G_0(t-x) f_i(t_{(i)}) dt.
\end{eqnarray*}
Note that for all $i=1, \ldots, |I|$
\begin{eqnarray*}
\int G_0(t-x) f_i(t_{(i)}) dt &=& 
\int \Bigl[\prod_{j\in I_i} g(t_j-x_j)\Bigr]
f_{i}(t_{(i)}) d t_{(i)} 
\\
\int G_{i,h}(t-x) f_i(t_{(i)}) dt &=& \int \Bigl[\prod_{j\in I_i}
\frac{1}{h_j} g\Bigl(\frac{t_j-x_j}{h_j}\Bigr)\Bigr] f_i(t_{(i)}) dt_{(i)}
\\
\int G_{j,h}(t-x) f_i(t_{(i)}) dt &=& \int \Bigl[\prod_{j\in I_i}
g(t_j-x_j)\Bigr] f_i(t_{(i)}) dt_{(i)},\;\;\;j\ne i.
\end{eqnarray*}
Combining these equalities
we obtain 
\[
\int K_{\theta}(t-x) F(t) dt = 
\int \Bigl[\prod_{j\in I_i}
\frac{1}{h_j} g\Bigl(\frac{t_j-x_j}{h_j}\Bigr)\Bigr] f_i(t_{(i)}) dt_{(i)},
\]
and
\[
\begin{array}{l}
{\displaystyle
B_{\theta}(x)
=
\sum_{i=1}^{|I|}
\int \Bigl[\prod_{j\in I_i}
\frac{1}{h_j} g\Bigl(\frac{t_j-x_j}{h_j}\Bigr)\Bigr] 
[f_i(t_{(i)})- f_i(x_{(i)})] dt_{(i)}
}
\\
{\displaystyle
=
\sum_{i=1}^{|I|}
\int \Bigl[\prod_{j\in I_i}
\frac{1}{h_j} g\Bigl(\frac{t_j-x_j}{h_j}\Bigr)\Bigr] 
\Bigl[f_i(t_{(i)})- f_i(x_{(i)})-
\sum_{s=1}^{l_i} \frac{1}{s!}
\sum_{|\blk|=s} D^{\blk} f_i(x_{(i)}) (t_{(i)} - x_{(i)})^{\blk}
\Bigr] dt_{(i)},
}
\end{array}
\]
where the last equality follows from 
the fact that
\[
\int \prod_{j\in I_i} \frac{1}{h_j} g\Bigl(\frac{t_j-x_j}{h_j}\Bigr)
(t_{(i)}-x_{(i)})^{\blk} dt_{(i)}=0,\;\;\;
\forall |\blk|: |\blk|=1,\ldots, l_i,\;\;i=1, \ldots, |I|,
\]
see Assumption~G(i). 
Because $f_i \in H_{{|I_i|}}(\beta_i, L_i)$, we obtain
\begin{eqnarray*}
|B_\theta (x)|
\leq  \sum_{i=1}^{|I|} L_i \int \Bigl|
\prod_{j\in I_i} \frac{1}{h_j} g\Bigl(\frac{t_j-x_j}{h_j}\Bigr)\Bigr|
 |t_{(i)}-x_{(i)}|^{\beta_i} dt_{(i)}
\leq
 \sum_{i=1}^{|I|} L_i \|g\|_1^{{|I_i|}} \sum_{j\in I_i} h_j^{\beta_i}.
\end{eqnarray*}
as claimed.
\epr
We quote the following result from \citeasnoun{talagrand}
that is repeatedly used in the proof of Lemma~\ref{lem:exp-inequalities} below.
\begin{lemma}\label{lem:talagrand}
Consider a centered Gaussian process $(X_t)_{t\in T}$.
Let $\sigma^2=\sup_{t\in T} EX_t^2$. Consider the intrinsic
semi--metric  
$\rho_X$ on $T$ given by $\rho_X^2 (s, t)=\rE(X_s-X_t)^2$.
Assume that for some constant $A>\sigma$, some $v>0$ and some
$0\leq \e_0\leq \sigma$ we have
\[
\e<\e_0 \;\;\Rightarrow\;\;N(T, \rho_X, \e)\leq \Bigl(\frac{A}{\e}\Bigr)^v,
\]
where $N(T, \rho_X, \e)$ is the smallest number of balls
of radius $\e$ needed to cover $T$. 
Then for $u\geq \sigma^2[(1+\sqrt{v})/\e_0]$
we have
\[
\rP\Bigl(\sup_{t\in T} X_t \geq u \Bigr) \leq 
\Bigr(\frac{KAu}{\sqrt{v}\sigma^2}\Bigr)^v \Phi
\Bigl(\frac{u}{\sigma}\Bigr),
\]
where $K$ is universal constant,  and $\Phi(u)=\frac{1}{\sqrt{2\pi}}\int_u^\infty 
e^{-s^2/2} ds$.
\end{lemma}
\begin{lemma}\label{lem:exp-inequalities}
Let Assumptions A, K0 and K2 hold. Then
for any $\kappa\geq 1+\sqrt{\frac{d+m}{\gamma}}$ one has
\begin{equation}\label{eq:exp-1}
\rP\big\{\sup_{\theta\in\Theta} \|\tilde{Z}_\theta (\cdot)\|_\infty \geq \kappa\big\}
\leq N [C_1 \bar{L} R\kappa]^{(d+m)/\gamma} \exp\{-\kappa^2/2\},
\end{equation}
where $C_1$ is an absolute constant.
\par
Furthermore, for any $\kappa\geq 1+\sqrt{\frac{d+2m}{\gamma}}$ one has
\begin{equation}\label{eq:exp-2}
\rP\big\{\sup_{(\theta, \nu)\in\Theta\times\Theta} 
\|\tilde{Z}_{\theta, \nu} (\cdot)\|_\infty \geq \kappa\big\}
\leq N^2 [C_2 M({\cal K})\bar{L} R\kappa]^{(d+2m)/\gamma} \exp\{-\kappa^2/2\},
\end{equation}
where $C_2$ is an absolute constant.
\end{lemma}
\pr  1$^0$. First we prove (\ref{eq:exp-1}). 
Recall our notation:
\[
Z_\theta (x) = \int K_\theta (t, x) W(dt),\;\;
\sigma_\theta (x)=\|K_\theta (\cdot, x)\|_2,\;\;
\tilde{Z}_\theta (x)= \sigma_{\theta}^{-1}(x) Z_{\theta} (x). 
\]
By Assumption~A, $\theta=(\theta_1, \theta_2)\in \Theta_1\times \Theta_2$.
Because the set $\Theta_1$ is finite, throughout the proof we keep
$\theta_1\in \Theta_1$ fixed.
For brevity, we will write
$\theta=(\theta_1, \theta_2)$, $\theta^\prime=(\theta_1, \theta_2^\prime)$,
$u=(x, \theta_2)$, $u^\prime=(x^\prime, \theta_2^\prime)$. Also with a slight abuse of notation we write
$Z(u)$, $\tilde{Z}(u)$ and $\sigma(u)$ for
$Z_{\theta}(x)$, $\tilde{Z}_\theta (x)$ and $\sigma_{\theta}(x)$ respectively.
The same notation with $u$ replaced by $u^\prime$
will be used for the corresponding quantities depending on $u^\prime$.
\par
Consider the random process 
$\{Z(u), u\in U\}$. Clearly, it has zero 
mean
and variance $\rE Z^2(u) = \sigma^2(u)$. 
Let $\rho_Z$ denote the intrinsic
semi--metric of $\{Z(u), u\in U\}$; then
\begin{eqnarray*}
\rho_Z (u, u^\prime) &:=& [\rE |Z(u) - Z(u^\prime)|^2 ]^{1/2}
\\
&=&
\|K_{(\theta_1, \theta_2)}(\cdot, x) -
K_{(\theta_1, \theta_2^\prime)}(\cdot, x^\prime)\|_2
\\
&\leq& 
\bar{L} |u-u^\prime|^\gamma,
\end{eqnarray*}
where the last inequality follows from Assumption~K2.
\par
Now consider the random process 
$\{\tilde{Z}(u), u\in U\}$. Let 
$\underline{\sigma}=\inf_{u\in U} \sigma(u)$;
then
\begin{eqnarray}
\rho_{\tilde{Z}} (u, u^\prime) &:=& 
\big[\rE|\tilde{Z}(u)-\tilde{Z}(u^\prime)|^2]^{1/2}
\nonumber
\\
&=&
\Bigl[
\rE\big| \frac{Z(u)}{\sigma(u)}   
- \frac{Z(u^\prime)}{\sigma(u^\prime)}\big|^2
\Bigr]^{1/2}
\nonumber
\\
&\leq&
\frac{1}{\sigma (u)} \rho_Z (u, u^\prime) + \sigma(u^\prime)
\big| \frac{1}{\sigma (u)}-\frac{1}{\sigma (u^\prime)}\big|
\nonumber
\\
&\leq& 
\underline{\sigma}^{-1} 
\big[ \rho_Z (u, u^\prime)+
| \sigma (u)-\sigma(u^\prime)|\big]
\nonumber
\\
&\leq & 2\underline{\sigma}^{-1} \rho_Z (u, u^\prime) 
\leq 2 (\mes\{\cD\})^{1/2} \bar{L} |u-u^\prime|^\gamma.
\label{eq:intrinsic-1}
\end{eqnarray}
Here we have taken into account that $\underline{\sigma}\geq (\mes \{\cD\})^{-1/2}$, and
\begin{eqnarray*}
|\sigma (u) - \sigma(u^\prime)| &=& |\, \|K_{\theta}(\cdot, x)\|_2 -\|K_{\theta^\prime}(\cdot, x^\prime)\|_2\,|
\\
&\leq & \|K_{\theta}(\cdot, x) -K_{\theta^\prime}(\cdot, x^\prime)\|_2
=\rho_Z (u, u^\prime).
\end{eqnarray*}
\par
It follows from (\ref{eq:intrinsic-1}) that the covering number
$N(U, \rho_{\tilde{Z}}, \eta)$ of the index set $U=\cD_0\times \Theta_2$
with respect to the intrinsic semi--metric $\rho_{\tilde{Z}}$
does not exceed $[c_1 \bar{L}R \eta^{-1}]^{(d+m)/\gamma}$, where
$c_1$ is an absolute constant. Then using 
the exponential inequality of Lemma~\ref{lem:talagrand}
[with $v=(d+m)/\gamma$, $A=c_1\bar{L}R$ and $\sigma=\e_0=1$], 
and summing over all $\theta_1\in\Theta_1$ we obtain
(\ref{eq:exp-1}).
\par
2$^0$. Now we turn to the proof of (\ref{eq:exp-2}). We recall that
\begin{eqnarray*}
Z_{\theta, \nu}(x)-Z_\nu(x) &=&
\int \big[K_{\theta, \nu}(t, x)- K_\nu(t, x)\big] W(dt),
\\
\sigma_{\theta, \nu}(x) &=& \|K_{\theta, \nu}(\cdot, x)-K_{\nu}(\cdot, x)\|_2,
\end{eqnarray*}
where $K_{\theta, \nu}(\cdot, \cdot)$ is defined in 
(\ref{eq:convoluted-kernel}).
We keep $\theta_1, \nu_1 \in \Theta_1$ fixed, and denote
$\theta=(\theta_1, \theta_2)$, 
$\theta^\prime=(\theta_1, \theta_2^\prime)$, $\nu=(\nu_1, \nu_2)$,
$\nu^\prime=(\nu_1, \nu_2^\prime)$. We also denote
$V=\cD_0\times \Theta_2\times \Theta_2$, $v=(\theta, \nu, x)$,
$v^\prime (\theta^\prime, \nu^\prime, x^\prime)$, and 
consider the Gaussian random processes
$\{\zeta (v), v\in V\}$ and $\{\tilde{\zeta}(v), v\in V\}$, where
\[
\zeta(v)= Z_{\theta, \nu}(x) - Z_{\nu}(x),\;\;\;
\tilde{\zeta}(v)= 
\tilde{\sigma}_{\theta, \nu}^{-1}(x) [Z_{\theta, \nu}(x) - Z_{\nu}(x)].
\]
Let $\rho_\zeta$ and $\rho_{\tilde{\zeta}}$ be the intrinsic semi--metrics
of these processes. Similarly to (\ref{eq:intrinsic-1}), it is straightforward
to show that $\rho_{\tilde{\zeta}}(v, v^\prime)\leq 2\rho_{\zeta}(v, v^\prime)$,
and our current goal is to bound $\rho_\zeta(v, v^\prime)$ from above.
\par
We have
\begin{eqnarray*}
\rho_\zeta (v, v^\prime) &=& \big[\rE |\zeta(v) -\zeta(v^\prime)|^2\big]^{1/2}
\\
&=&
\|K_{\theta, \nu}(\cdot, x)-K_\nu(\cdot, x) - 
K_{\theta^\prime, \nu^\prime}(\cdot, x^\prime)+K_{\nu^\prime}(\cdot, x^\prime)
\|_2
\\
&\leq &
\|K_\nu(\cdot, x) - K_{\nu^\prime}(\cdot, x^\prime)\|_2+
\|K_{\theta, \nu}(\cdot, x)- K_{\theta^\prime, \nu^\prime}(\cdot, x^\prime)\|_2
= J_1+J_2.
\end{eqnarray*}
By Assumption~K2
\[
J_1\leq \bar{L} |v-v^\prime|^\gamma.
\]
Let $\widehat{g}(\cdot, x)$ be the Fourier transform of a function 
$g: \cD\times \cD_0 \to \rR^1$ with respect to the first argument, i.e.,
\[
\widehat{g}(\omega, x)= \int g(t, x) \exp\{2\pi i \omega^T t\} dt,\;\;
\forall x\in \cD_0.
\]
Then, by construction,
$\widehat{K}_{\theta, \nu}(\cdot, x)= \widehat{K}_\theta (\cdot, x)
\widehat{K}_\nu (\cdot, x)$, and 
\begin{eqnarray*}
J_2 &=& \| \widehat{K}_\theta(\cdot, x) \widehat{K}_{\nu}(\cdot, x)-
\widehat{K}_{\theta^\prime} (\cdot, x^\prime) 
\widehat{K}_{\nu^\prime}(\cdot, x^\prime)\|_2
\\
&\leq& \| \,[\widehat{K}_\theta (\cdot, x)- \widehat{K}_{\theta^\prime}
(\cdot,x^\prime)]
\widehat{K}_\nu (\cdot, x)\|_2 +
\| \,[\widehat{K}_\nu (\cdot, x)- \widehat{K}_{\nu^\prime}(\cdot,x^\prime)]
\widehat{K}_{\theta^\prime} (\cdot, x^\prime)\|_2
\\
&\leq &
\| K_{\nu}(\cdot, x)\|_1 \, \|K_\theta (\cdot, x)- K_{\theta^\prime}
(\cdot,x^\prime)\|_2
+
\|K_{\theta^\prime} (\cdot, x^\prime)\|_1\, 
\|K_\nu (\cdot, x)- K_{\nu^\prime}(\cdot,x^\prime)\|_2
\\
&\leq & 2M({\cal K}) \bar{L} |u-u^\prime |^\gamma,
\end{eqnarray*}
where we have used Assumptions K0 and K2. Combining upper bounds for $J_1$ and $J_2$ we get
$\rho_\zeta (v, v^\prime) \leq [1+2M({\cal K})] \bar{L} |v-v^\prime|^\gamma$,
and finally
\begin{equation}\label{eq:intrinsic-2}
\rho_{\tilde{\zeta}}(v, v^\prime) \leq 2[1+2M({\cal K})] \bar{L} |v-v^\prime|^\gamma.
\end{equation}
\par
It follows from (\ref{eq:intrinsic-2}) that the covering number
$N(V, \rho_{\tilde{\zeta}}, \eta)$ of the index set 
$V=\cD_0\times \Theta_2\times\Theta_2$
with respect to the intrinsic semi--metric $\rho_{\tilde{\zeta}}$
does not exceed $[c_2 M({\cal K})\bar{L}R \eta^{-1}]^{(d+2m)/\gamma}$, where
$c_2$ is an absolute constant. Then noting that 
$\sup_v {\rm var}(\tilde{\zeta}(v))\leq 1$,   using 
the exponential inequality of Lemma~\ref{lem:talagrand}
[with $v=(d+2m)/\gamma$, $A=c_2M({\cal K})\bar{L}R$ and $\sigma=\e_0=1$], 
and summing over all $(\theta_1, \nu_1)\in\Theta_1\times \Theta_1$ we obtain~(\ref{eq:exp-2}).

\bibliographystyle{agsm}

\end{document}